\documentclass[final,a4paper]{jmsj}
\usepackage{xcolor}
\usepackage{mathrsfs,amsmath,amssymb,amsthm}
\usepackage[dvipdfmx]{graphicx}
\usepackage{ulem}
\usepackage[]{showlabels}

\newcommand{\R}{\mathbb R}
\newcommand{\dd}{d}

 \newcommand{\Sv}{S}
\newcommand{\Tv}{R}

\newcommand{\TT}{T_t}

\newcommand{\bnormal}{\boldsymbol{\nu}}
\newcommand{\bfxi}{\boldsymbol{\xi}}
\newcommand{\bfal}{\boldsymbol{\alpha}}
\newcommand{\bfs}{\mathbf{s}}

\newcommand{\bfa}{\mathbf{a}}
\newcommand{\bfb}{\mathbf{b}}

\newtheorem{theorem}{Theorem}
\newtheorem{lemma}[theorem]{Lemma}
\newtheorem{definition}[theorem]{Definition}
\newtheorem{remark}[theorem]{Remark}

\begin{document}

\title[Liouvilles's formulae and Hadamard variation]
{Liouville's formulae and Hadamard variation with respect to general domain perturbations}

\author{Takashi \textsc{Suzuki}}
\address{Center of Mathematical Modeling and Data Science, \\
Osaka University, \\
Toyonaka, 560-8531, Japan}
\email{suzuki@sigmath.es.osaka-u.ac.jp}

\author{Takuya \textsc{Tsuchiya}}
\address{Graduate School of Science and Engineering, \\
Ehime University, \\
Matsuyama 790-8577, Japan}
\email{tsuchiya.takuya.plateau@kyudai.jp}

\subjclass[2020]{Primary 35J25; Secondary 35R35}

\keywords{the Green function, domain perturbations, the Hadamard
variation, Liouville's formulae, the Neumann boundary condition}

\recdate{20XX}{11}{99}
\revdate{20XX}{1}{99}

\begin{abstract}
We study Hadamard variations with respect to general domain perturbations, particularly for the Neumann boundary condition. They are derived from new  Liouville's formulae  concerning the transformation of volume and area integrals.  Then, relations to several geometric quantities are discussed; differential forms and the second fundamental form on the boundary. 
\end{abstract}

\maketitle

\section{Introduction} 

Our purpose is to establish Liouville's formulae on volume and area
integrals and derive Hadamard variations with respect to general domain
perturbations for the Neumann boundary condition. 

Let $\Omega$ be a bounded domain in $d$-dimensional
Euclidean space $\R^d$, $d \ge 2$.  If its boundary $\partial\Omega$
is represented as graphs of Lipschitz functions, $\Omega$ is called a
\textit{Lipschitz domain}.  For the detailed definition of Lipschitz
domains, see \cite{SuzTsu16} and the references there in. If $\Omega$
is a Lipschitz domain, the smoothness of $\partial\Omega$ is
denoted by $C^{0,1}$.
We suppose that $\partial\Omega$ is divided into two non-overlapped closed 
sets $\gamma^0$ and $\gamma^1$ satisfying 
\begin{align}
   \gamma^0 \cup \gamma^1= \partial\Omega,
   \quad
   \gamma^0 \cap \gamma^1 = \emptyset.
  \label{assum-gamma}
\end{align}
This assumption yields that $\gamma^i$, $i=0,1$, do not have their
boundaries in $\partial\Omega$.  

We study the Poisson problem with
the mixed boundary condition:
\begin{align}
  - \Delta z = f \ \mbox{in $\Omega$}, \quad
      z = \varphi \ \mbox{on $\gamma^0$}, \quad
      \frac{\partial z}{\partial \bnormal} = \psi \ \mbox{on $\gamma^1$},
   \label{eqn:poisson}
\end{align}
where
$\Delta=\sum_{i=1}^d \frac{\partial^2}{\partial x_i^2}$ stands for the
Laplacian and $\bnormal$ is the unit outer normal vector on $\partial\Omega$. The standard theory of  elliptic PDE tells us that if $\partial\Omega$ is $C^2$, or $C^{1,1}$ more weakly, which means that it is Lipschitz continuous up to the first derivatives, and $f \in L^2(\Omega)$, $\varphi \in H^2(\Omega)$, and $\psi \in H^1(\Omega)$, equation \eqref{eqn:poisson} admits a
unique solution $z \in H^2(\Omega)$. (See Section \ref{seca1} below for the case of Lipschitz domains.) 

This solution admits the representation  
\begin{equation}
   z(y) = \int_{\Omega} N(x,y) f(x) \ \dd x
    - \int_{\gamma^0} \varphi(x)\frac{\partial N}{\partial \bnormal_x}
   (x,y) \ \dd s_x + \int_{\gamma^1}N(x,y)\psi(x) \ \dd s_x 
   \label{solution-ex}
\end{equation}
for $y\in \Omega$, where $ds_x$ is the surface element and $N(x,y)$ is the Green's function. Thus, given $y\in \Omega$, if $\Gamma(x)$ denotes the fundamental solution of $\Delta$: 
\begin{equation}
  \Gamma(x) =\left\{ 
    \begin{array}{ll} 
     - \frac1{2\pi}\log \left| x\right|, & d = 2, \\
      \frac1{(d-2)\omega_d}|x|^{2-d}, & d \ge 3,
    \end{array} \right.
  \label{4.0}
\end{equation}
where $\omega_d$ is the volume of the unit ball in $\R^d$, and if $u=u(x)$ is the solution to 
\begin{equation}
   \Delta u = 0 \ \mbox{in $\Omega$}, \quad
   u= - \Gamma(\cdot -y) \ \mbox{on $\gamma^0$}, \quad  
   \frac{\partial u}{\partial \bnormal}
   = - \frac{\partial \;}{\partial \bnormal}
   \Gamma(\cdot - y) \ \mbox{on $\gamma^1$}, 
    \label{4}
\end{equation}
this $N(x,y)$ is given by  
\begin{equation}
  N(x,y)= \Gamma(x-y) + u(x).
   \label{5}
\end{equation}

We take a family of domain perturbations parametrized
by $t$, $|t| \ll 1$, which is denoted by $\{\TT\}$ (the exact definition will be given in
Section~\ref{sec:DomainPertur} below).  Suppose that the boundaries $\gamma^i$, $i=0,1$, are
mapped onto $\gamma_t^i$, $i = 0,1$, respectively, by $\TT$: 
\[ \TT(\gamma^i)=\gamma^i_t, \quad i=0,1. \] 
Then the Green's function on $\Omega_t$ is defined by 
\begin{equation}
N(x,y,t) = \Gamma(x - y) + u(x,t),
 \label{4.1}
\end{equation}
where $u=u(x,t)$ is the solution to  
\begin{equation}
  \Delta u(\cdot,t) = 0 \ \mbox{in $\Omega_t$}, \quad
  u(\cdot,t) = - \Gamma(\cdot-y) \ \mbox{on $\gamma_t^0$}, \quad
  \frac{\partial u}{\partial \nu}(\cdot,t)
  = - \frac{\partial}{\partial \nu}
  \Gamma(\cdot - y) \ \mbox{on $\gamma_t^1$}.
   \label{4.5}
\end{equation}
By this definition, it holds that $\Omega_0=\Omega$, $N(x,y,0)=N(x,y)$, and $u(x,0)=u(x)$.

Given $x,y\in\Omega$, we have $x,y\in \Omega_t$ 
for $\vert t\vert \ll 1$. Then the Hadamard variation of the Green's function $N(x,y)$ is defined by
\begin{equation}
   \delta N(x,y)= \frac{\partial N}{\partial t}
     (x,y,t)\Big|_{t=0} 
    = \frac{\partial u}{\partial t}
     (x,t)\Big|_{t=0}, \quad
    (x,y)\in \Omega \times \Omega.
 \label{6} 
\end{equation}
The second variation is defined similarly: 
\begin{equation}  
\delta^2N(x,y)= \frac{\partial^2N}{\partial t^2}(x,y,t)\Big|_{t=0} =\frac{\partial^2u}{\partial t^2}(x,t)\Big|_{t=0}, \quad (x,y)\in \Omega \times \Omega. 
 \label{77}
\end{equation} 
These notions are classical, but to clarify the meaning of these derivatives, including their existence, is one of our main aims. The other is to prescribe a class of domains which admits these limits. The Lipschitz domain is a main target, from the viewpoint of numerical computations in engineering, which are often concerned on polygons on the plane. Hence we are taking the applications to free boundary problems \cite{st05} or shape optimizations \cite{azegami} in mind. These problems in engineering induce the third motivation of ours, study on the general domain perturbation of the domain. Thus we continue our previous work \cite{SuzTsu16} on the Dirichlet boundary condition. 

So far, the domain perturbation has been often introduced by the deformation of $\partial \Omega$ as in 
\begin{equation} 
\partial \Omega_t: x+t\cdot \delta \rho(x)\nu_x, \ x\in \partial \Omega, 
 \label{nperturbation}
\end{equation} 
where $\delta \rho(x)$ is a smooth function of $x\in \partial\Omega$. This method, which may be called the normal perturbation, does not always work for the general Lipschitz domain, for example, if $\partial\Omega$ has a corner. The dynamical perturbation introduced by \cite{SuzTsu16} may fit more the Lipschitz domain. It is given by $\TT x=X(t)$ for the solution $X=X(t)$ to 
\begin{equation} 
\frac{dX}{dt}=v(X), \quad \left. X\right\vert_{t=0}=x\in \overline{\Omega}, 
 \label{7}
\end{equation} 
where $v=v(x)$ is a Lipschitz continuous vector field defined on an open neighbourhood of $\overline{\Omega}$. 

In \cite{SuzTsu16} we have studied the first and the second Hadamard variations under the general perturbation of Lipschitz domains for the case $\gamma^1=\emptyset$, that is, for the Dirichlet boundary condition. This paper is devoted to the general case of $\gamma^1\neq \emptyset$. Hence it is concerned on the Neumann boundary condition essentially, and extends the classical result of \cite{gs52} on (\ref{nperturbation}) for $d=2,3$. Meanwhile we carefully examine the regularity of the domain admitting these variations. 

Our strategy is a systematic use of Liouville's formulae concerning the variation of volume and area integrals under the tranformation of variables. Actually, we have derived Liouville's first volume and area formulae in the general form to treat the Dirichlet boundary condition in \cite{SuzTsu16}.  Here we formulate the second  formulae of these integrals to study the Neumann boundary condition. These formulate are concerned on the second derivatives and are to be used for numerical computations adapting Newton's method. From the viewpoint of pure mathematics, on the other hand, a role of the second fundamental form of the boundary in the second Hadamard variation is clarified for general domain perturbations. The other topic is the derivation of the first and the second area formulae via the transformation of differential forms. 

This paper is composed of four sections and two appendices. Taking preliminaries in $\S$\ref{sec2}, we show Liouville's volume and area formulae in $\S$\ref{sec3}. Then $\S$\ref{sec4} is devoted to the Hadamard variation with respect to general domain perturbations. The first appendix, $\S$\ref{seca2}, is devoted to the derivation of Liouville's area formulae via differential forms.  In the second appendix, $\S$\ref{seca3}, we show a form of Liouville's second area formula represented by the second fundamental form of $\partial \Omega$. 

 Several formulae on Hadamard variations of eigenvalues noticed by \cite{gs52}, such as the harmonic concavity of the first eigenvalue, will be generalized in our forthcoming paper, with rigorous proof of the existence of these variations.

\section{Preliminaries}\label{sec2} 

\subsection{Poisson equation on Lipschitz domains}\label{seca1}

The fundamental property of the Lipschitz domain is the following fact \cite[Theorem 3, p.127]{eg92}. 
\begin{theorem}\label{thm00}
If $\Omega\subset \R^d$ is a Lipschitz domain, then the set of functions $C^\infty(\overline{\Omega})$ is dense in $W^{1,p}(\Omega)$ for $1\leq p<\infty$, where 
\begin{equation} 
C^\infty(\overline{\Omega})=\{ v:\overline{\Omega}\rightarrow \R \mid \exists \tilde v\in C_0^\infty(\R^d), \ \left. \tilde v\right\vert_{\overline{\Omega}}=v\}. 
 \label{not1}
\end{equation} 
\end{theorem} 

This theorem ensures the validity of the trace operator to
$\partial\Omega$ (\cite[Theorem 1, p.133]{eg92}), and then the unique
solvability of (\ref{eqn:poisson}) for the Lipschitz domain $\Omega$
holds as in the case of $\gamma^1=\emptyset$. Hence we have the
following theorem similarly to \cite{SuzTsu16}. Note that the spaces
$C^{0,1}(\partial\Omega)$ and $H^\theta(\partial\Omega)$ for
$0<\theta<1$ are well-defined by the local chart because
$\partial\Omega$ is Lipschitz continuous. 

\begin{theorem}\label{thm0}
If $\Omega \subset \R^d$ is a bounded Lipschitz domain satisfying (\ref{assum-gamma}), there arise the following facts: First, the trace operator $v\in C^\infty(\overline{\Omega})\mapsto \left. v\right\vert_{\partial \Omega}\in C^{0,1}(\partial\Omega)$ is extended as 
\[  
v\in H^1(\Omega) \ \mapsto \ \left. v\right\vert_{\partial\Omega}\in H^{1/2}(\partial\Omega).  
\] 
Then there arise the isomorphisms 
\[  
v\in H^1(\Omega)/H^1_0(\Omega) \ \mapsto \ \left. v\right\vert_{\partial\Omega}\in H^{1/2}(\partial\Omega)  \] 
and   
\[  
v\in V/H^1_0(\Omega) \ \mapsto \ \left. v\right\vert_{\gamma^1}\in H^{1/2}(\gamma^1)  \] 
for 
$V=\{ v\in H^1(\Omega) \mid \left. v\right\vert_{\gamma^0}=0\}$. Second, the normal derivative of $v\in H^1(\Omega)$ on $\partial\Omega$ is defined in the sense of  
\[
   \frac{\partial v}{\partial \bnormal}\in H^{-1/2}(\partial\Omega)\,
   = H^{1/2}(\partial\Omega)',
\] 
provided that $\Delta v\in H^1(\Omega)'$, 
where $\Delta$ is taken in the sense of distributions in $\Omega$.
Hence there arises that 
\[
   \left\langle \varphi, \frac{\partial v}{\partial \nu}
  \right\rangle_{H^{1/2}(\partial\Omega),
   H^{-1/2}(\partial\Omega)} = (\nabla v, \nabla\varphi)_{L^2(\Omega)}
  +\langle \varphi, \Delta v\rangle_{H^1(\Omega), H^{1}(\Omega)'} 
\] 
for any $\varphi\in H^1_0(\Omega)$. Similarly, if $\Delta v\in V'$, the normal derivative of $v\in V$ on $\gamma^1$ is defined in the sense of  
\[ \frac{\partial v}{\partial \bnormal}\in H^{-1/2}(\gamma^1), \] 
and it holds that   
\[ \left\langle \varphi, \frac{\partial v}{\partial \nu}\right\rangle_{H^{1/2}(\gamma^1), H^{-1/2}(\gamma^1)}=(\nabla v, \nabla\varphi)_{L^2(\Omega)}+\langle \varphi, \Delta v\rangle_{V, V'} \] 
for any $\varphi\in V$. Finally, given $f\in V'$, $\varphi\in H^{1/2}(\gamma^0)$, and $\psi\in H^{-1/2}(\gamma^1)$, there is a unique solution $z\in H^1(\Omega)$ to (\ref{eqn:poisson}). Hence this $z$ satisfies  
\[ \left. z\right\vert_{\gamma^0}=\varphi \] 
and 
\[
 (\nabla z, \nabla \zeta)=(f,\zeta)+\langle \zeta, \psi
 \rangle_{V,V'}, \quad \forall \zeta \in V.
 \] 
\end{theorem} 

The Green's function $N(x,y)$ of (\ref{eqn:poisson}) is now defined by (\ref{4})-(\ref{5}). It satisfies 
\[ N(\cdot,y)\in V, \ -\Delta N(\cdot, y)=\delta(\cdot-y), \ \left. \frac{\partial N}{\partial \bnormal}(\cdot,y)\right\vert_{\gamma^1}=0, \quad \forall y\in \Omega, \] 
where $\delta(x)$ stands for the delta function and $\frac{\partial N}{\partial \bnormal}(\cdot,y)$ on $\gamma^1$ is taken as an element in $H^{-1/2}(\gamma^1)$. Hence the solution $z\in H^1(\Omega)$ to (\ref{eqn:poisson}) for $\varphi\in H^{1/2}(\gamma^0)$ and $\psi \in H^{-1/2}(\gamma^1)$ admits the representation    
\begin{equation} 
   z(y) = (N(\cdot,y), f) - \left\langle \varphi, \frac{\partial N}{\partial \bnormal}
   (\cdot,y) \right\rangle_{\gamma^0} + \langle N(\cdot,y), \psi \rangle_{\gamma^1}, \quad \forall y\in \Omega. 
 \label{8}
\end{equation} 
Here and henceforth, $( \ , \ )$ and $\langle \ , \ \rangle_{\gamma^i}$, $i=0,1$, denote the inner product in $L^2(\Omega)$ and the paring between $H^{1/2}(\gamma^i)$ and $H^{-1/2}(\gamma^i)$, respectively. 

We note that the above $H^1$ theory to (\ref{eqn:poisson}) is valid even if $\sigma\equiv \gamma^0\cap \gamma^1\neq \emptyset$ in (\ref{assum-gamma}), provided that the $(d-1)$ dimensional Hausdorff measure of $\sigma$ vanishes.

\subsection{Differentiations on $\partial\Omega$}

We continue to suppose that $\Omega \subset \R^d$ is a bounded Lipschitz domain satisfying (\ref{assum-gamma}). It follows from Rademacher's
theorem that the tangent space $T_x(\partial\Omega)$ and 
the unit outer normal vector $\bnormal$ exist for almost every $x\in\partial\Omega$.  At such $x \in \partial\Omega$, we take the orthonormal moving frame
\begin{align*}
  \{\bfs_1, \cdots, \bfs_{d-1}, \bnormal\}
\end{align*}
with positive orientation, where 
$\{\bfs_1, \cdots, \bfs_{d-1}\}$ is an orthonormal
frame (with positive orientation) of the tangent space
$T_x(\partial\Omega)$.  Equalities concerning the derivatives of Lipschitz functions below are valid almost everywhere, although it is not mentioned each time.  

Let $\gamma=\gamma^i$, $i=0,1$, and $\nabla=\nabla_x$. If $\gamma$ is $C^{1,1}$, the above $\bfs_1, \cdots, \bfs_{d-1}, \bnormal$ are $C^{0,1}(\gamma)$. If $\gamma$ is $C^{2,1}$, which means that it is Lipschitz continuous up to the derivatives of the second order, these vectors $\bfs_1, \cdots, \bfs_{d-1}, \bnormal$ are in $C^{1,1}(\gamma)$. In this case, if $f$ is a $C^{1,1}$ function in a neighbourhood of $\gamma$, then we obtain (\cite[Lemma~10, Corollary~11]{SuzTsu16}), 
\begin{eqnarray}
   \nabla ^2f & = & \sum_{i=1}^{d-1}(\nabla \bfs_{i})^T
 \frac{\partial f}{\partial s_{i}} + 
 (\nabla\bnormal)^T\frac{\partial f}{\partial\nu}
 + \sum_{i,j=1}^{d-1}\bfs_{i}\otimes \bfs_{j}
  \frac{\partial^{2}f}{\partial s_{i}\partial s_{j}} \nonumber\\
  & & + \sum_{i=1}^{d-1}
  (\bfs_{i}\otimes \bnormal + \bnormal\otimes \bfs_{i})
  \frac{\partial^{2}f}{\partial s_{i}\partial \nu}
 + \bnormal\otimes \bnormal \frac{\partial^{2}f}{\partial \nu^{2}}  
 \label{eqn:hessian}, 
\end{eqnarray} 
and 
\begin{equation} 
\Delta f = (\nabla\cdot \bnormal)\frac{\partial f}{\partial \nu}
 + \frac{\partial^{2}f}{\partial \nu^{2}} + 
 \sum_{i=1}^{d-1}\frac{\partial^2f}{\partial s_i^2}, \quad
   \nabla\cdot\bnormal = \sum_{i=1}^{d-1}\kappa_i 
 \label{eqn:laplacian}
\end{equation}
on $\Gamma$, where $\nabla^2 f$ is the Hesse matrix of $f$, 
\begin{eqnarray}
   (\nabla\bnormal)^T & 
  = & \sum_{i=1}^{d-1}\kappa_i \bfs_i\otimes \bfs_i
    + \sum_{i=1}^{d-1}\sum_{j\neq i}
   \left(\frac{\partial \bnormal}{\partial s_i}
   \cdot \bfs_j\right) \bfs_j\otimes \bfs_i,  
   \label{nablanu} \nonumber\\
  (\nabla \bfs_j)^T &
  = & - \kappa_j \bnormal\otimes \bfs_j - \sum_{i \neq j}
   \left(\frac{\partial \bnormal}{\partial s_i}
   \cdot \bfs_j\right) \bnormal \otimes \bfs_i,
   \label{nablas}
\end{eqnarray}
and $\kappa_i$ is the sectional curvature of $\gamma$ along $\bfs_j$,
$j = 1, \cdots, d-1$.

In the general case of the bounded Lipschitz domain $\Omega$, each 
$\gamma^i$, $i=0,1$, forms a Lipschitz manifold without boundaries.
Then the Stokes theorem ensures  
\begin{equation}
   \int_{\gamma^i} \dd \omega = 0, 
   \quad i = 1, 2,
    \label{StokesThm}
\end{equation}
where $\omega$ is a Lipschitz continuous differential form of order $d-2$ and $\dd\omega$
is its exterior derivative.

Let $F$ and $g$ be Lipschitz continuous functions on $\gamma^1$, and $H$ be $C^{1,1}$ in a neighbourbood of $\gamma^1$. Let, furthermore, $\omega_1(F,H)$ and $\omega_2(F,g,H)$ be differential forms with
order $d-2$ defined by
\begin{eqnarray*}
& & \omega_1(F,H) = \sum_{i=1}^{d-1} (-1)^{i+1} F H
      \, \dd \bfs_1 \wedge \cdots \wedge
    \widehat{\dd s_i} \wedge \cdots \wedge \dd \bfs_{d-1} \\
& & \omega_2(F,g,H) = \sum_{i=1}^{d-1} (-1)^{i+1} F g
    \frac{\partial H}{\partial \bfs_i} 
    \, \dd \bfs_1 \wedge \cdots \wedge
    \widehat{\dd s_i} \wedge \cdots \wedge \dd \bfs_{d-1},
\end{eqnarray*}
where $\widehat{\dd \bfs_i}$ means excluding of $\dd \bfs_i$ and
$\frac{\partial \;}{\partial \bfs_i}$ is the directional derivative
along $\bfs_i$.  Using \eqref{StokesThm}, we obtain 
\[ 
  \langle \dd \omega_1(F,H), \, 1 \rangle_{\gamma^1} = 
  \langle \dd \omega_2(F,g,H), \, 1 \rangle_{\gamma^1} = 0,
\] 
and hence the following lemma.  Here, the tangential gradient $\nabla_\tau$ on $\partial\Omega$ is
defined by
\begin{equation}
   \nabla_\tau= \sum_{i=1}^{d-1} \frac{\partial \;}{\partial \bfs_i},
   \label{tan-grad}
\end{equation}
which is independent of the choice on the
orthonormal coordinate $\{\bfs_1, \cdots, \bfs_{d-1}\}$. 

\begin{lemma}
It holds that
\begin{eqnarray}
   \sum_{i=1}^{N-1} \left\langle \frac{\partial F}{\partial \bfs_i},
   \,  H \right\rangle_{\gamma^1} & = & - \sum_{i=1}^{N-1} \left\langle F, \, 
    \frac{\partial H}{\partial \bfs_i}\right\rangle_{\gamma^1},
    \label{dif-form0} \\ 
  \langle \nabla_{\tau} F, \,g \nabla_{\tau} H
   \rangle_{\gamma^1} & = &  - \left\langle F, \, \sum_{i=1}^{N-1}
    \frac{\partial \;}{\partial \bfs_i} \left(
     g \frac{\partial H}{\partial \bfs_i}
      \right) \right\rangle_{\gamma^1}.
    \label{dif-form}
\end{eqnarray} 
\end{lemma}

\subsection{Domain perturbations}\label{sec:DomainPertur}

Given the bounded Lipschitz domain $\Omega\subset {\R}^d$, 
let 
\begin{equation} 
\TT : \Omega\rightarrow \Omega_t= \TT(\Omega), \ |t| \ll 1, \quad 
 {T}_0=I
 \label{13}
\end{equation} 
be a family of bi-Lipschitz homeomorphisms. 

\begin{definition} 
The family of deformations $\{\TT\}$ of $\Omega$ in (\ref{13}) is said to be differentiable 
if $\TT x$ is continuously differentiable in $t$ for every $x\in\Omega$ and the mappings 
\[ 
  \frac{\partial}{\partial t} D\TT, \ 
  \frac{\partial}{\partial t} (D\TT)^{-1}: 
  \Omega \rightarrow \R^d 
\]  
are uniformly bounded, where $D\TT$ denotes the Jacobi matrix of $\TT:\Omega\rightarrow \Omega_t$. It is said to be twice differentiable if $\TT x$ is continuously differentiable twice in $t$ for every $x\in \Omega$ and the mappings
\[ 
   \frac{\partial^2}{\partial t^2}D\TT, \
   \frac{\partial^2}{\partial t^2}(D\TT)^{-1}: 
  \Omega\rightarrow \R^d 
  \] 
are uniformly bounded.
\end{definition}

The dynamical perturbation (\ref{7}) is once and twice differentiable if $v\in C^{0,1}(\tilde \Omega; \R^d)$ and $v\in C^{1,1}(\tilde \Omega; \R^d)$, respectively, where $\tilde \Omega$ is an open neighbourhood of $\overline{\Omega}$. If $\{\TT\}$ is, say, twice differentiable, the vector fields $\Sv$ and $\Tv$ defined by 
\begin{equation}
   \Sv = \frac{\partial \TT}{\partial t}\bigg|_{t=0}, \quad
   \Tv = \frac{\partial^2 \TT}{\partial t^2}\bigg|_{t=0}
 \label{taylor0}
\end{equation} 
are Lipschitz continuous on $\overline{\Omega}$. Then the family $\{\TT\}$ admits the Taylor expansion in $C^{0,1}(\overline{\Omega})$, 
\begin{equation}
  \TT =I + t \Sv + \frac{1}{2} t^2 \Tv + o(t^2),  \quad t\rightarrow 0,   
  \label{T-taylor}
\end{equation}
if it is twice differentiable, where $I$ denotes the identity mapping. Then we put  
\begin{equation}
 \delta\rho = \frac{\partial \TT}{\partial t} \bigg|_{t=0}
 \cdot  \bnormal = \Sv\cdot \bnormal, \quad 
 \delta^{2}\rho= \frac{\partial^{2} \TT}{\partial t^{2}}
   \bigg|_{t=0} \cdot \bnormal = \Tv \cdot \bnormal, 
    \label{14} 
\end{equation}
recalling that $\bnormal$ is the unit outer normal vector on $\partial\Omega$. 

For the normal perturbation, it holds that 
\begin{equation} 
S=\rho\bnormal \ \mbox{on $\partial\Omega$}, \quad R=0,  
 \label{consistent}
\end{equation} 
and therefore, $\delta\rho$ in (\ref{14}) is consistent to that in (\ref{nperturbation}). If $\{\TT\}$ is a dynamical perturbation defined by (\ref{7}), then we obtain 
\begin{equation} 
S=v, \quad R=(v\cdot \nabla)v 
 \label{19}
\end{equation} 
and hence  
\[ \delta \rho=v\cdot \bnormal, \quad \delta^2\rho=[(v\cdot \nabla)v]\cdot \bnormal. \] 

Careful treatments of the domain perturbation are necessary if $\gamma^0\cap \gamma^1\neq \emptyset$ is admitted in (\ref{assum-gamma}), such as the non-pealing-off condition used in \cite{st11}. This case is left in a future study. 

\subsection{Jacobi matrix and its derivatives}

Given the $d \times d$ matrices $A=(a_{ij})$ and $B=(b_{ij})$,
their inner product $A:B$ and the associated Frobenius norm $\|A\|_F$
are defined by
\begin{align*}
A:B= \sum_{i,j=1}^d a_{ij}b_{ij}, \quad
\|A\|_F^2 = A:A = \sum_{i,j=1}^d a_{ij}^2. 
\end{align*}
The Jacobi matrix $D\TT$ of the bi-Lipschitz 
homeomorphism $\TT$ is defined by 
\begin{align*}
D\TT= \left(\frac{\partial \TT^i}{\partial x_j} 
\right)_{i,j=1, \cdots, d}, 
\end{align*}
where 
\[ 
\TT= \left(\TT^1, \cdots, \TT^d\right)^T.
\] 

We continue to suppose that $\Omega \subset \R^d$ is a bounded Lipschitz domain and $\{\TT\}$, $\vert t\vert \ll 1$, is a family of bi-Lipschitz deformations of $\Omega$. 

\begin{lemma}\label{deri-det0}
It holds that  
\begin{equation}
\lim_{t \to 0} D\TT = I, \quad \lim_{t \to 0} \frac{\partial\;}{\partial t} (D\TT) = D\Sv, \quad  
\lim_{t\to 0} \frac{\partial^2}{\partial t^2}(DT_t)=DR
 \label{15}
\end{equation} 
uniformly on $\Omega$ if $\{ \TT\}$ is differentiable, and 
\begin{equation} 
\lim_{t \to 0} \frac{\partial\;}{\partial t}(D\TT)^{-1}
= - D \Sv, \quad \lim_{t \to 0}\frac{\partial^2\;}{\partial t^2} (D\TT)^{-1} 
= 2(D\Sv)^2 - D\Tv 
  \label{lim-t}
\end{equation}
uniformly on $\Omega$ if $\{\TT\}$ is twice differentiable, furthermore. Here, $\Sv$ and $\Tv$ are vector fields on $\Omega$ defined by \eqref{T-taylor}, and $D\Sv$ and $D\Tv$ are the Jacobi matrices of $\Sv$ and $\Tv$, respectively.
\end{lemma}

\begin{proof}
Convergences below are uniform on $\Omega$. The limits in (\ref{15}) are obvious by \eqref{T-taylor}.  For (\ref{lim-t}), we differentiate 
\begin{equation} 
I = (D \TT) (D \TT) ^ {-1} 
 \label{17}
\end{equation} 
once and twice with respect to $t$, to obtain
\begin{eqnarray*}
  O & = & \left(\frac{\partial \;}{\partial t}
   (D\TT)\right) (D \TT)^{-1}
  + (D \TT)\frac{\partial \;}{\partial t} (D \TT)^{-1}, \\
  O & = & \left(\frac{\partial^2 \;}{\partial t^2} (D\TT)\right)
   (D \TT)^{-1}
+ 2 \frac{\partial \;}{\partial t} (D\TT)
\frac{\partial \;}{\partial t} (D \TT)^{-1}
+ (D \TT)\frac{\partial^2 \;}{\partial t^2} (D \TT)^{-1}. 
\end{eqnarray*}
With $t \to 0$, there arises that 
\[ 
 O = D\Sv + \left. \frac{\partial \;}{\partial t} (D \TT)^{-1}
\right|_{t = 0}, \quad 
 O = D\Tv  + 2(D\Sv)(-D\Sv) 
 + \left. \frac{\partial^2 \;}{\partial t^2} (D \TT)^{-1}
\right|_{t = 0} 
\] 
by 
\[ \lim_{t \to 0} (D\TT)^{-1} = I \] 
derived from (\ref{17}). Then (\ref{lim-t}) follows. 
\end{proof}

The following lemma is given in \cite[Theorem~6]{SuzTsu16}.
\begin{lemma}\label{deri-det1}
It holds that 
\begin{equation}
\label{eqno325}
  \left. \frac{\partial}{\partial t} \, \mathrm{det} \, D\TT
   \right\vert_{t=0} = \nabla\cdot \Sv 
\end{equation}
uniformly on $\Omega$ if $\{ \TT\}$ is differentiable, and  
\begin{equation}
\label{eqno268}
  \left. \frac{\partial^{2}}{\partial t^{2}} \, \mathrm{det} \,
   D\TT  \right\vert_{t=0}=\nabla\cdot \Tv 
  + (\nabla\cdot \Sv)^2 - (D\Sv)^T:D\Sv
\end{equation}
 uniformly on $\Omega$ if $\{ \TT\}$ is twice differentiable. 
 \label{thm:jacobi}
\end{lemma}

\section{Liouville's formulae}\label{sec3} 

\subsection{First formulae} 

This section is devoted to several Liouville's formulae concerning volume and area integrals under general perturbation of Lipchitz domains. They are used to derive Hadamard variations associated with the Neumann boundary condition in the following section. The formulae given below are concerned on general domain perturbations. They are new, and have their own interests. 

We continue to suppose that $\Omega \subset \R^d$ is a bounded Lipschitz domain and $\{\TT\}$, $\vert t\vert <\varepsilon$, is a  family of deformations of $\Omega$. We write $\langle \ , \ \rangle_{\partial \Omega}$ either for the paring between $H^{1/2}(\partial\Omega)$ and $H^{-1/2}(\partial\Omega)$, or, for the inner product in $L^2(\partial\Omega)$. Differentiations in $t$ of the volume and area integrals below are taken in the classical sense, unless otherwise stated.  

The first volume formula follows from (\ref{eqno325}) and a transformation of variables as in \cite[Theorem~1]{SuzTsu16}. Note that $\delta\rho$ defined by (\ref{14}) is Lipschitz continuous on $\partial\Omega$.  Let $Q\subset \R^{d+1}$ be the non-cylindrical domain defined by 
\begin{equation} 
Q=\bigcup_{\vert t\vert<\varepsilon}\Omega_t\times \{ t\}. 
 \label{21}
\end{equation} 

\begin{theorem}[first volume formula]
\label{liouville-first}
If $\Omega \subset \R^d$ is a bounded Lipschitz domain, $\{ \TT\}$ is differentiable, $c\in C^{0,1}(\overline{Q})$, and $c_t$ is continuous on $\overline{Q}$, it holds that 
\[ 
  \frac{\dd \;}{\dd t}\left. 
  \int_{\Omega_t}c \ \dd x \right|_{t=0}
   = \int_{\Omega}\dot c \ \dd x +
   \langle  c_0, \, \delta\rho \rangle_{\partial\Omega},  
\] 
where $c_0=c(\cdot,0)$ and $\dot c=c_t(\cdot,0)$. 
\end{theorem}

\begin{remark}
\upshape
Without the last requirement of $c_t\in C(\overline{Q})$, we have 
\begin{equation} 
  \frac{\dd I}{\dd t} 
   = \int_{\Omega_t}c_t \ \dd x +\int_{\partial\Omega_t}c \
 \left(\frac{\partial\TT}{\partial t}\cdot \nu_t\right) \  \dd s_t 
 \label{nshka}
\end{equation}
in the sense of distributions in $t$, where $\dd s_t$ denotes the area element of $\partial\Omega_t$ and 
\begin{equation} 
I=\int_{\Omega_t}c \ dx. 
 \label{ai}
\end{equation} 
In other words, the mapping $t\mapsto I(t)$ is locally absolutely continuous, and equality (\ref{nshka}) holds for almost every $t$. Then the above condition $c_t\in C(\overline{Q})$ is used to take the initial trace at $t=0$ in (\ref{nshka}).  
\end{remark} 

\begin{remark}
\upshape
In \cite{SuzTsu16} the above $c=c(x,t)$ is required to be extended outside $\overline{Q}$. This extension is always possible within the category of Lipschitz continuous functions because $Q\subset \R^{d+1}$ is a  Lipschitz domain by the dynamical extension described later. Hence if $c\in C^{0,1}(\overline{Q})$, there is $\tilde c\in C^{0,1}(\hat Q)$ such that $\left. \tilde c\right\vert_{Q}=c$, where $\hat Q$ is an open set containing $\overline{Q}$. 
\end{remark} 

The following form with less regularity of $c=c(x,t)$ is used later for computations of the Hadamard variation. Note that Theorem \ref{liouville-first} is a direct consequence of this theorem. 

\begin{theorem}\label{thm8}  
Let $\Omega\subset \R^d$ be a bounded Lipschitz domain and $\{ \TT\}$ be differentiable. Given $c\in L^{1}(Q)$ with $\nabla c\in L^1(Q; \R^d)$, put 
\begin{equation} 
b(x,t)=c(\TT x,t), \quad \beta(x,t)=\nabla c(\TT x,t)\cdot\frac{\partial \TT}{\partial t}(x), 
 \label{32.0}
\end{equation} 
and assume $b\in C^1(-\varepsilon, \varepsilon; W^{1,1}(\Omega))$. 
Then it holds that 
\[ \alpha(x,t)\equiv c_t(\TT x, t)\in L^1(\tilde Q), \quad \tilde Q=\Omega\times (-\varepsilon, \varepsilon). \] 
Assume, furthermore, $\beta \in C^0(-\varepsilon, \varepsilon; L^1(\Omega))$. Then it holds that $\alpha\in C^0(-\varepsilon, \varepsilon; L^1(\Omega))$, $I\in C^1(-\varepsilon, \varepsilon)$ for $I=I(t)$ defined by (\ref{ai}), and 
\begin{equation} 
\left. \frac{dI}{dt}\right\vert_{t=0} = \int_\Omega \dot c+\nabla\cdot (c_0S) \ dx, 
 \label{conclusion}
\end{equation} 
where $c_0=c(\cdot, 0)$, $\dot c=c_t(\cdot,0)$, 
\end{theorem} 
\begin{proof} 
Since $c\in L^1(Q)$, $\nabla c\in L^1(Q; \R^d)$, and $b\in C^{1}(-\varepsilon, \varepsilon; W^{1,1}(\Omega))$, we obtain  $\beta \in L^1(\tilde Q)$ and 
\begin{equation} 
b_t=\alpha+\beta \quad \mbox{in $\tilde Q=\Omega\times (-\varepsilon, \varepsilon)$}  
 \label{35x}
\end{equation}   
in the sense of distributions, and hence $\alpha \in L^1(\tilde Q)$. If $\beta\in C^0(-\varepsilon, \varepsilon; L^1(\Omega))$, furthermore, it holds that $\alpha \in C^0(-\varepsilon, \varepsilon; L^1(\Omega))$, and therefore, the volume integral on the right-hand of (\ref{conclusion}) is definite by 
\[ c_t(\TT x, t), \ \nabla c(\TT x,t)\cdot \frac{\partial \TT}{\partial t}(x) \ \in C^0(-\varepsilon, \varepsilon; L^1(\Omega)). \]   
 
Since  
\begin{equation} 
\int_{\Omega_t}c \ \dd x=\int_\Omega b(x,t) \ \mbox{det} \ D\TT (x) \ \dd x  
 \label{35}
\end{equation} 
and $\{ \TT\}$ is differentiable, we have 
\begin{eqnarray*} 
\frac{dI}{dt} & = & \int_{\Omega} b_t+b\frac{\partial}{\partial t}\det D\TT \ dx  = \int_\Omega \alpha+\beta+b\frac{\partial}{\partial t}\det D\TT \ dx \\ 
& = & \int_\Omega c_t(\TT x, t)+\nabla c(\TT x,t)\cdot \frac{\partial \TT}{\partial t}(x)+c(\TT x, t)\frac{\partial}{\partial t}\det D\TT x \ dx  
\end{eqnarray*}  
with the continuity of the right-hand side in $t$. Then the result follows with $t=0$,  
\[ \left.\frac{dI}{dt}\right\vert_{t=0}=\int_\Omega \dot c+(\nabla c_0)\cdot S+c_0\nabla\cdot S \ dx, \] 
from (\ref{eqno325}). 
\end{proof} 

\begin{remark}
\upshape
The differentiabilities of $b$ and $c$ in $t$ stand for those of the Lagrange and the Euler ones, respectively. Later in the study of Hadamard variations, we use the fact that the differentiability of the former implies that of the latter by the above theorem.    
\end{remark}

To formulate Liouville's area formula, we define a $d$-dimensional Lipschitz-manifold in $\R^{d+1}$ by   
\begin{equation} 
\Gamma=\bigcup_{\vert t\vert<\varepsilon}\partial \Omega_t\times \{t\}. 
 \label{gamma}
\end{equation}
Here we employ the method of dynamical extension for the proof of the following lemma. 

\begin{lemma}\label{remark10} 
Each $f\in C^{0,1}(\Gamma)$ is extended as an element in $C^{0,1}(\tilde \Gamma)$,  where $\tilde \Gamma$ is an open neighbourhood of $\Gamma'$ in $\R^{d+1}$ for 
\begin{equation} 
\Gamma'=\bigcup_{\vert t\vert <\varepsilon'}\partial\Omega_t\times \{t\}, \quad 0<\varepsilon'<\varepsilon. 
 \label{27}
\end{equation} 
\end{lemma} 

\begin{proof} 
We take a smooth vector field $v=v(z)$, $z=(x,t)$, in $\R^{d+1}$. 
Then, we let $Z=Z(z, s)$ to be the solution to 
\[ \frac{dZ}{ds}=v(Z), \quad \left. Z\right\vert_{s=0}=z\in \Gamma. \] 
There is $0<\delta_0\ll 1$, such that 
\[ \tilde \Gamma=\{ Z(z,s) \mid z\in \Gamma, \ \vert s\vert <\delta_0\} \] 
forms an open neighbourhood of $\Gamma'$. 
Furthermore, by the uniqueness of the solution of the
ordinary differential equation, the orbits, 
 $\{ {\cal O}_z\mid  z\in \Gamma\}$ with
${\mathcal O}_z=\{ Z(z,s) \mid \vert s \vert <\delta_0\}$, 
do not intersect each other, 
and form a tubular neighourhood of $\Gamma$.
 Given $f\in C^{0,1}(\Gamma)$, now we put 
\[ \tilde f(Z)=f(z), \quad Z=Z(z,s)\in \tilde \Gamma, \]  
to obtain $\tilde f\in C^{0,1}(\tilde \Gamma)$.  
\end{proof} 

\begin{remark} \upshape
Given $\bfa\in C^{0,1}(\Gamma; \R^d)$, we thus obtain its extension, denoted by the same symbol, $\bfa \in C^{0,1}(\tilde \Gamma; \R^d)$, which assures $\bfa_t \in L^\infty(\tilde\Gamma; \R^d)$. This property, however, does not imply $\bfa_t(\cdot,t) \in L^\infty(\partial\Omega_t; \R^d)$, $\vert t\vert< \varepsilon$, because of the discrepancy of the dimensions of $\tilde \Gamma$ and $\partial\Omega_t$, the boundary of $\Omega_t$. Note that $\nabla\cdot\bfa(\cdot,t)\in L^\infty(\partial\Omega_t)$ holds for each $t$ because $\bfa \in C^{0,1}(\Gamma; \R^d)$ implies $\bfa(\cdot, t)\in C^{0,1}(\partial\Omega_t; \R^d)$. Hence we require for both $\bnormal \cdot \bfa_t$ and $\nabla\cdot \bfa$ to be continuous on $\Gamma$, besides $\bfa \in C^{0,1}(\Gamma)$, in Lemma \ref{vector-first} below.  
\end{remark} 

Recalling that $ds_t$ denotes the area element of $\partial\Omega_t$, we put $ds=ds_0$. The outer unit normal vector $\bnormal=\bnormal(\cdot,t)$ is defined almost everywhere on $\partial \Omega_t$ for any $t$, because $\Omega_t$ is a Lipschitz domain. 

\begin{lemma}\label{vector-first}
If $\Omega \subset \R^d$ is a bounded Lipschitz domain, $\{ \TT\}$ is differentiable, $\bfa\in C^{0,1}(\Gamma; \R^d)$, and both $\bnormal\cdot\bfa_t$ and $\nabla\cdot \bfa$ are continuous on $\Gamma$, it holds that 
\[ 
 \frac{\dd \;}{\dd t}\left.
  \int_{\partial\Omega_t} \bnormal\cdot \bfa \ \dd s_t \right|_{t=0}
  = \int_{\partial\Omega} (\bnormal\cdot \bfa_t)(\cdot,0)\, \dd s
   +  \langle \nabla \cdot \bfa(\cdot,0), 
       \delta \rho \rangle_{\partial\Omega}.
\] 
\end{lemma}
\begin{proof}
This lemma is reduced to an identity, valid to any $\bfa \in W^{1,1}(\Gamma)$, that is, 
\begin{equation} 
 \frac{\dd \;}{\dd t}
  \int_{\partial\Omega_t} \bnormal\cdot \bfa \ \dd s_t 
  =  \int_{\partial\Omega_t} \biggl(
   \bnormal\cdot \bfa_t 
   +  (\nabla \cdot \bfa) \frac{\partial \TT}{\partial t}\cdot \bnormal
    \biggr)  \  \dd s_t,
 \label{24}
\end{equation}
in the sense of distributions in $t$. Hence in (\ref{24}), the boundary integral of the left-hand side is locally absolutely continuous in $t\in (-\varepsilon, \varepsilon)$, and equality holds for almost all $t$.  Here, the required regularity ensures the Lipschitz continuity of $\bnormal\cdot \bfa_t$, $\nabla\cdot \bfa$, and $\frac{\partial \TT}{\partial t}\cdot \bnormal$ on $\Gamma$, and also the differentiability in $t$ of the bi-Lipschitz diffeomorphism $\TT: \partial\Omega\rightarrow \partial \Omega_t$. We thus obtain the continuity of the right-hand side of (\ref{24}) in $t$, and hence the conclusion. 

By Lemma \ref{remark10} there is an extension of $\bfa$, denoted by the same symbol, such that $\bfa \in C^{0,1}(\tilde \Gamma)$, Then $C^\infty(\overline{\tilde \Gamma})$ is dense in $W^{1,1}(\tilde \Gamma)$ by Theorem \ref{thm00}. We may assume, therefore, $\bfa\in C^\infty(\overline{\Gamma})$ to verify (\ref{24}), recalling the notion (\ref{not1}). In this case this differentiation in $t$ is to be valid in the classical sense, and the equality is to hold for all $t$. Below we describe the proof of this fact just for $t=0$ to make the description simple. Hence we show the lemma for $\bfa=\bfa(x,t)$ smooth in $\tilde \Gamma$, the open  neighbourhood of $\Gamma'$ in $\R^{d+1}$ defined by (\ref{27}).  

For this purpose, we extend $\bfa$ to a smooth vector field $\tilde \bfa$ on $\overline{Q'}$ for 
\[ Q'=\bigcup_{\vert t\vert<\varepsilon'}\Omega\times \{t\}, \] 
that is, $\tilde \bfa=\bfa \varphi$, where $0 \leq \varphi=\varphi(x,t)\leq 1$ is a smooth function supported in $\tilde \Gamma$ and is equal to $1$ near $\Gamma'$. Thus the proof of this lemma is reduced for smooth $\bfa$, say $\bfa\in C^{1,1}(\overline{Q'}, \R^d)$. 

Now we use the transformation of variables, 
\begin{equation} 
\bfb(x,t) = \bfa(y,t), \quad y=\TT x. 
 \label{trv}
\end{equation} 
It holds that $D_x \bfb(x,t) = D_y \bfa(y, t) D\TT(x)$ and hence 
\[ 
 D_y \bfa(y,t) = D_x \bfb(x, t) (D\TT(x))^{-1},
\]  
where $D_x \bfb$ and $D_y\bfa$ denote the Jacobi matrices of $\bfb$ and $\bfa$ with respect to
$x$ and $y$, respectively. Then Green's formula implies 
\begin{eqnarray}
\int_{\partial\Omega_t} (\bnormal\cdot\bfa)(y,t) \ \dd s_t  & = &  
 \int_{\Omega_t} \nabla_y\cdot\bfa(y,t) \ \dd y 
  = \int_{\Omega_t} \mathrm{tr}\left[
   D_y \bfa(y,t)\right] \dd y \nonumber\\
  & = & \int_\Omega \mathrm{tr}\left[
    D_x \bfb(x, t) (D_x\TT(x))^{-1} \right]
    \mathrm{det}(D\TT)(x) \ \dd x, 
   \label{25}
\end{eqnarray}
where $\mathrm{tr}\, X$ denotes the trace of the matrix $X$.  Note that the right-hand side of (\ref{25}) is differentiable in $t$ by the above reduction of $\bfa\in C^{1,1}(\overline{Q}, \R^d)$. 

In fact, we have 
\begin{eqnarray}  
& & \frac{\partial}{\partial t}\left\{\mathrm{tr}\left[
    D_x \bfb(\cdot, t) (D_x\TT)^{-1} \right]
    \mathrm{det}(D\TT)\right\} \nonumber\\ 
& & \quad =\mathrm{tr}\left[
    D_x \bfb(\cdot, t) (D_x\TT)^{-1} \right]_t \mathrm{det}(D\TT)+ \mathrm{tr}\left[
    D_x \bfb(\cdot, t) (D_x\TT)^{-1} \right]
    (\mathrm{det}(D\TT))_t.   
   \label{22}
\end{eqnarray}  
Then Lemmas \ref{deri-det0} and \ref{deri-det1} imply  
\[ 
D_x \bfb(x,t) \to D_x\bfa(x ,0), \quad (\mathrm{det}(D\TT))_t(x) \to (\nabla\cdot\Sv)(x) 
\] 
and 
\begin{eqnarray*} 
& & D_x \bfb_t(x,t)(D\TT(x))^{-1} + D_x \bfb(x,t)
  \frac{\partial\;}{\partial t} (D\TT(x))^{-1} \\ 
  & & \quad \to \ D_x \bfb_t(x,t)\bigm|_{t=0} - D_x\bfa(x,0)(D\Sv)(x) 
\end{eqnarray*} 
as $t \to 0$, uniformly in $x\in \overline{\Omega}$. Here, since 
\begin{eqnarray*}
& & D_x \bfb_t(x,t) \\ 
& & \quad = D_x\left( D_y \bfa(\TT x,t)
  \frac{\partial \TT}{\partial t}(x) + \bfa_t(\TT x,t)\right) \\ 
& & \quad = \left(D_y^2 \bfa(\TT x,t) D\TT(x)\right) 
     \frac{\partial \TT}{\partial t}(x) + D_y \bfa(\TT x,t) \frac{\partial \;}{\partial t} (D\TT(x))
  + D_y \bfa_t(\TT x,t) D\TT(x)
\end{eqnarray*} 
it holds that  
\begin{equation} 
D_x \bfb_t(x,t)\bigm|_{t=0}=D_x^2 \bfa(x,0)\Sv(x) + D_x \bfa(x,0) D\Sv(x) + D_x \bfa_t(x,0), 
 \label{limit_bt}
\end{equation} 
where $D_x^2 \bfa$ is the third-order tensor consisting of the second
derivative of the components of $\bfa$. Hence there arises that  
\[ D_x^2\bfa(\cdot,0)\Sv=\left( \sum_k\frac{\partial^2a^i}{\partial x_j\partial x_k}(\cdot,0)S^k \right)_{i,j=1, \cdots, d} \] 
for $\bfa=(a^i)$ and $\Sv=(S^i)$. 

Gathering these observations, we obtain
\begin{eqnarray*}
& & \frac{\dd \;}{\dd t}\int_{\partial\Omega_t} 
    (\bnormal\cdot\bfa)(\cdot,t) \ \dd s_t  \biggr|_{t=0} \\ 
& & \quad = \frac{\dd \;}{\dd t}\left. 
     \int_\Omega \mathrm{tr}\left[
    D_x \bfb(x, t) (D_x(\TT x))^{-1} \right]
    \mathrm{det} \ D\TT (x) \ \dd x \right|_{t=0} \\  
& & \quad = \int_\Omega \mathrm{tr}[D_x^2\bfa(\cdot,0)\Sv + D_x\bfa(\cdot,0)DS
+ D_x\bfa_t(\cdot,0)] + \mathrm{tr}[D_x\bfa(\cdot,0)(-DS)] \\ 
& & \qquad \qquad \qquad + \mathrm{tr}[D_x\bfa(\cdot,0)] (\nabla\cdot\Sv) \ \dd x \\
& & \quad = \int_\Omega \mathrm{tr}[D_x^2\bfa(\cdot,0)\Sv] + (\nabla_x\cdot\bfa)(\cdot,0) (\nabla\cdot\Sv)+ (\nabla_x\cdot\bfa_t)(\cdot,0)  \ \dd x \\ 
& & \quad = \int_\Omega \nabla_x\cdot [(\nabla_x \cdot \bfa) (\cdot,0)) \Sv] + (\nabla_x \cdot \bfa_t)(\cdot,0) \ \dd x \\
& & \quad = \int_{\partial\Omega} (\bnormal\cdot \bfa_t)(\cdot,0) 
   + (\nabla\cdot\bfa)(\cdot,0) (\Sv \cdot \bnormal) \ \dd s, 
\end{eqnarray*}
and hence the conclusion. 
\end{proof}

\begin{remark}\upshape 
Lemma \ref{vector-first} ensures an extension of \cite[Lemma~13]{SuzTsu16}. Here, it is not necessary to assume $c=\nabla\cdot\bfa$ with a vector field $\bfa$, $C^2$ in a neighbourhood of $\Gamma$. The outer unit normal vector $\bnormal$ in the following theorem is in 
\[ \bnormal\in C^{0,1}(\Gamma, S^{d-1}), \quad S^{d-1}=\{ \zeta\in \R^d \ \mid \vert \zeta\vert=1\} \] 
from the assumption. Then we obtain $\bnormal(\cdot,t)\in C^{0,1}(\partial\Omega_t, S^{d-1})$ and hence $(\nabla\cdot \bnormal)(\cdot,t)\in L^\infty(\partial \Omega_t)$ for each $t$, which is equal to the mean curvature of $\partial \Omega_t$.
\end{remark} 

\begin{theorem}[first area formula]
\label{liouville-first-are}
If $\Omega\subset \R^d$ is a bounded $C^{1,1}$ domain, $\{ \TT\}$ is differentiable, $c\in C^{0,1}(\Gamma)$, and both $c_t$ and $\frac{\partial c}{\partial \bnormal}$ are continuous on $\Gamma$, it holds that
\[ 
   \frac{\dd \;}{\dd t}\left.
   \int_{\partial\Omega_t}c \ \dd s_t \right|_{t=0}
    = \int_{\partial\Omega}\dot c \ \dd s + 
    \left\langle (\nabla\cdot\bnormal)c_0 +\frac{\partial c_0}{\partial
 \bnormal},  \delta\rho \right\rangle_{\partial\Omega},
\] 
where $c_0=c(\cdot,0)$ and $\dot c=c_t(\cdot, 0)$. 
\end{theorem}
\begin{proof} 
Lemma \ref{remark10} ensures  an extension of $\bnormal\in C^{0,1}(\Gamma; S^{d-1})$ to 
\[ \bnormal \in C^{0,1}(\tilde \Gamma; S^{d-1}), \] 
where $\tilde \Gamma$ is an open neighbourhood of $\Gamma'$ in $\R^{d+1}$. Then we apply Lemma \ref{vector-first} to $\bfa= \bnormal c\in C^{0,1}(\Gamma)$.  

In fact, it holds that   
\begin{equation} 
 \bnormal \cdot \bnormal_t = 0 \quad \mbox{on $\Gamma$}  
 \label{by}
\end{equation} 
by $\vert \bnormal\vert^2 = 1$ in $\tilde\Gamma$, and therefore, $\bnormal\cdot\bfa_t= c_t$ is continuous on $\Gamma$ from the assumption. Then we obtain the result by  
\[ \nabla \cdot \bfa = (\nabla\cdot\bnormal)c + \bnormal\cdot\nabla c \quad \mbox{on $\Gamma$}.  \] 
\end{proof}

\subsection{Second formulae} 

We turn to Liouville's second formulae. In the following formula on the volume integral, the last term vanishes for the dynamical perturbation $\{ \TT\}$ defined by (\ref{7}), because of   
\[ \delta^2\rho-[(\Sv\cdot\nabla)\Sv ]\cdot\bnormal=(R-[\Sv\cdot\nabla \Sv])\cdot\bnormal \] 
and (\ref{19}). 

\begin{theorem}[second volume formula]
\label{second-liouville-thm}
If $\Omega\subset \R^d$ is a bounded Lipschitz domain, $\{ \TT\}$ is twice differentiable, $c\in C^{1,1}(\overline{Q})$, and $c_{tt}$ is continuous on $\overline{Q}$, it holds that 
\[ 
\frac{\dd^2 \;}{\dd t^2} \left. 
  \int_{\Omega_t}c \ \dd x \right|_{t=0} 
  = \int_\Omega \ddot c \  \dd x
    + \langle 2\dot c + \nabla\cdot c_0\Sv, \delta\rho
    \rangle_{\partial\Omega} + \langle c_0, \delta^2\rho -[(\Sv\cdot\nabla)\Sv]\cdot\bnormal
   \rangle_{\partial\Omega}, 
\] 
where $c_0=c(\cdot,0)$, $\dot c=c_t(\cdot,0)$, and $\ddot c=c_{tt}(\cdot,0)$. 
\end{theorem}

\begin{proof}
The proof is reduced to the case of $c\in C^{\infty}(\overline{Q})$ as in Lemma \ref{vector-first}. Then we differentiate the right-hand side of
\begin{equation*}
  \int_{\Omega_t}c \ \dd x =
  \int_{\Omega} c(\TT x,t)\, \mathrm{det}(D\TT(x)) \ \dd x 
\end{equation*}
twice in $t$. It follows that 
\begin{eqnarray*}
& & \frac{d^2}{dt^2} \int_{\Omega_t}c \ \dd x = \int_{\Omega} \Big[c_{tt} (\TT x,t)
    + 2 \nabla_x c_t(\TT x,t)
     \cdot \frac{\partial \;}{\partial t}\TT x +
     [\nabla^2_xc(\TT x,t)] \left(
    \left(\frac{\partial \;}{\partial t}\TT x\right)^2 \right) \\ 
& & \quad + 
    \nabla_xc(\TT x,t)\cdot
     \frac{\partial^2 \;}{\partial t^2}\TT x\Big]\mbox{det} (D\TT(x))+ c(\TT x,t)
    \frac{\partial^2 \;}{\partial t^2} \mathrm{det}D\TT(x) \\
& & \quad + 2 \Big[c_t(\TT x,t) + \nabla_x c(\TT x,t)\cdot
    \frac{\partial \;}{\partial t}\TT x\Big]
    \frac{\partial \;}{\partial t} \mathrm{det}D\TT (x) \ \dd x.
\end{eqnarray*}
Letting $t\to 0$, then we obtain
\begin{eqnarray*}
& & \left.\frac{d^2}{dt^2} \int_{\Omega_t}c \ \dd x\right\vert_{t=0}= 
   \int_\Omega \left[
    \ddot c + 2 \nabla\dot c \cdot\Sv +
    \nabla c_0\cdot \Tv
   + (\nabla^2 c_0)(\Sv^2) \right]  \\ 
& & \quad + 2
    (\dot c + \nabla c_0\cdot \Sv) 
    (\nabla\cdot\Sv) 
    + c_0 (\nabla\cdot \Tv
    + (\nabla\cdot\Sv)^2 - D\Sv^T:D\Sv) \ \dd x 
\end{eqnarray*}
by Lemma~\ref{deri-det1}.  Since the divergence formula implies 
\begin{eqnarray*}
& & \int_\Omega \nabla c_0\cdot\Tv
     + c_0 (\nabla \cdot \Tv) \ \dd x = 
   \langle c_0, \Tv \cdot\bnormal \rangle_{\partial\Omega} = 
   \langle c_0, \delta^2\rho\rangle_{\partial\Omega} \\
& & \int_\Omega \nabla \dot c\cdot\Sv
      + \dot c (\nabla \cdot \Sv) \ \dd x = 
  \langle \dot c, \Sv \cdot\bnormal \rangle_{\partial\Omega}=
  \langle \dot c, \delta \rho\rangle_{\partial\Omega}, 
\end{eqnarray*}
it follows that 
\[ 
\left. \frac{\dd^2 \;}{\dd t^2} 
   \int_{\Omega_t}c \ \dd x \right|_{t=0}  =
   \int_\Omega \ddot c \ \dd x
  + 2 \langle \dot c, \delta\rho \rangle_{\partial\Omega}
  + \langle c_0, \delta^2\rho \rangle_{\partial\Omega} + X + Y + Z 
\] 
for  
\begin{eqnarray*} 
& & X= \int_\Omega [\nabla^2c_0](\Sv^2) \dd x \\ 
& & Y= 2 \int_\Omega ((\nabla c_0)\cdot \Sv) (\nabla\cdot\Sv) \ \dd x \\ 
& & Z= \int_\Omega c_0((\nabla\cdot\Sv)^2 - D\Sv^T:D\Sv) \ \dd x.
\end{eqnarray*}
We simplify the terms $X$, $Y$, and  $Z$ furthermore. 

Let $\Sv = (S^1, \cdots, S^N)^T$ and 
\[ c_j=\frac{\partial c}{\partial x_j}, \quad c_{ij}=\frac{\partial^2c}{\partial x_i\partial x_j}, \quad  S^i_j=\frac{\partial S^i}{\partial x_j},  \] 
for simplicitly. First, the divergence formula implies  
\begin{eqnarray*}
 X + Y & = & \sum_{i,j} \int_\Omega c_{ij}S^iS^j + 2c_{i}S^iS_{j}^j \ dx \\ 
  & = & \langle \nabla c \cdot \Sv, \Sv\cdot\bnormal
     \rangle_{\partial\Omega} +
   \sum_{i,j} \int_\Omega -c_{i} S^i_{j} S^j + c_{i}S^iS_{j}^j \ dx 
\end{eqnarray*}
Second, if $S\in C^{1,1}(\overline{Q})$ we obtain 
\begin{eqnarray*}
  Z & = & \sum_{i,j} \int_\Omega c (S_{i}^iS_{j}^j -
   S_{i}^jS_{j}^i) \ dx  \\
& = & \sum_{i,j}\int_{\partial \Omega} c(S^iS^j_j- S^jS^i_j)\nu^i \ ds + \sum_{i,j}\int_\Omega -(cS^j_j)_iS^i+(cS^i_j)_iS^j \ dx \\ 
& = & \langle c(\nabla \cdot \Sv), \Sv\cdot\bnormal
        \rangle_{\partial\Omega} - \langle c, \bnormal \cdot [(\Sv \cdot \nabla)\Sv] 
        \rangle_{\partial\Omega}
      +\sum_{i,j}\int_\Omega -c_{i}S^i S_{j}^j+ c_{i} S^j S_{j}^i \ dx, 
\end{eqnarray*}
for $\bnormal=(\nu^1, \cdots, \nu^d)^T$, similarly, by  
\[ \sum_{i,j}\int_\Omega cS^iS^j_{ij} \ dx=\sum_{i,j}\int_\Omega cS^jS^i_{ij} \ dx. \]  
Hence it holds that 
\[ Z=\langle c(\nabla \cdot \Sv), \Sv\cdot\bnormal \rangle_{\partial\Omega} - \langle c, \bnormal \cdot [(\Sv \cdot \nabla)\Sv] \rangle_{\partial\Omega}+\sum_{i,j}\int_\Omega -c_{i}S^i S_{j}^j + c_{i} S^j S_{j}^i \ dx \]  
for general $S\in C^{0,1}(\overline{Q})$. 

Adding these equations, we have
\begin{eqnarray*}
  X + Y + Z & = & \langle \nabla c \cdot \Sv, \Sv\cdot\bnormal
     \rangle_{\partial\Omega} +
    \langle c (\nabla \cdot \Sv), \Sv\cdot\bnormal
        \rangle_{\partial\Omega}
   - \langle c, \, \bnormal \cdot (\Sv \cdot \nabla)\Sv
        \rangle_{\partial\Omega}  \\
 & = & \langle \nabla \cdot (c \, \Sv), \Sv\cdot\bnormal
     \rangle_{\partial\Omega}
   - \langle c, \, \bnormal \cdot (\Sv \cdot \nabla)\Sv
        \rangle_{\partial\Omega}.
\end{eqnarray*}
Gathering all equations above, the proof is complete.
\end{proof}

The following form with less regularity of $c$ is applicable to the Hadamard variation as in Theorem \ref{thm8}.   The proof is the same and is omitted.  

\begin{theorem}\label{lem12}
Let $\Omega\subset \R^d$ be a bounded Lipschitz domain and $\{ \TT\}$ be twice differentiable. Given $c=c(x,t)\in L^{1}(Q)$ with $\nabla c\in L^1(Q; \R^d)$, suppose 
\[ b\in C^2(-\varepsilon, \varepsilon; W^{1,1}(\Omega)), \quad \beta\in C^0(-\varepsilon, \varepsilon; L^1(\Omega)) \] 
for $b=b(x,t)$ and $\beta=\beta(x,t)$ defined by (\ref{32.0}).  Then it holds that $\alpha\in C^0(-\varepsilon, \varepsilon; L^1(\Omega))$, $I\in C^2(-\varepsilon, \varepsilon)$, and 
\[ \frac{d^2}{dt^2}I=\int_{\Omega}b_{tt} \ \mbox{det} \ D\TT +b \ \frac{\partial^2}{\partial t^2}\mbox{det} \ D\TT +2(\alpha+\beta) \ \frac{\partial}{\partial t}\mbox{det} \ DT_t \ dx
\] 
for $I=I(t)$ defined by (\ref{ai}). 
\end{theorem} 

Liouville's second area formula is derived from the
following lemma.

\begin{lemma}\label{vector-second}
If $\Omega$ is $C^{1,1}$, $\{\TT\}$ is twice differentiable, $\bfa\in C^{1,1}(\Gamma; \R^d)$, and both $\bnormal\cdot\bfa_{tt}$ and $\nabla\cdot\bfa_t$ are continuous on $\Gamma$, it holds that 
\begin{eqnarray*}
& & \frac{\dd^2 \;}{\dd t^2}\left. 
  \int_{\partial\Omega_t} \bnormal\cdot \bfa \ \dd s_t\right|_{t=0}
  = \int_{\partial\Omega} \bnormal\cdot \bfa_{tt}(\cdot,0) \ \dd s + \langle 2\nabla\cdot \bfa_t(\cdot,0)+ \nabla\cdot\left[(\nabla\cdot\bfa(\cdot,0))\Sv
   \right], \delta \rho \rangle_{\partial\Omega} \\ 
& & \quad + \langle \nabla\cdot\bfa(\cdot,0), (\Tv - (\Sv\cdot\nabla)\Sv)\cdot\bnormal\rangle_{\partial\Omega}.
\end{eqnarray*}
\end{lemma}
\begin{proof}
Again, the proof is reduced to the case of $\bfa \in C^\infty(\Gamma; \R^d)$. Then we use the transformation (\ref{trv}), 
\[ \bfb(x,t)= \bfa(y,t), \quad y=\TT x,  \]  
to reach (\ref{25}): 
\[ 
\int_{\partial\Omega_t} \bnormal\cdot\bfa(y,t) \dd s_y 
   = \int_\Omega \mathrm{tr}\left[
    D_x \bfb(x, t) (D_x\TT(x))^{-1} \right]
    \mathrm{det}(D\TT) \dd x.
\] 

Differentiating the right-hand side twice, here we obtain 
\begin{eqnarray*} 
& & \frac{\partial^2\;}{\partial t^2} \ \mathrm{tr}\bigl[
    D_x \bfb(x, t) (D_x\TT(x))^{-1} \bigr]
    \mathrm{det}(D\TT(x)) = \mathrm{tr}\left[D_x\bfb_{tt}(x,t) (D\TT(x))^{-1}\right]
    \mathrm{det}(D\TT(x)) \\ 
& & \quad + \mathrm{tr}\left[D_x\bfb(x,t)
  [(D\TT(x))^{-1}]_{tt}\right]
    \mathrm{det}(D\TT(x)) + \mathrm{tr}\left[D_x\bfb(x,t)
  (D\TT)^{-1}\right] [ \mathrm{det}(D\TT(x))]_{tt} \\ 
  & & \quad + 2\mathrm{tr}\left[D_x\bfb_t(x,t)[(D\TT(x))^{-1}]_t\right]
    \mathrm{det}(D\TT) + \mathrm{tr}\left[D_x\bfb_t(x,t)
  (D\TT(x))^{-1}\right] [\mathrm{det}(D\TT(x))]_t \\ 
  & & \quad + 2 \mathrm{tr}\left[D_x\bfb(x,t)
  [(D\TT(x))^{-1}]_t\right]
    [\mathrm{det}(D\TT(x))]_t 
 \end{eqnarray*} 
Since 
\[ D_x \bfb(x,t) = D_y \bfa(\TT x, t) (D\TT(x)), \] 
there arises that 
\begin{eqnarray*}
& &  D_x \bfb_{tt}(x,t) = 
   \left(D_y^3 \bfa(\TT x,t) D\TT(x)\right) 
     \left(\frac{\partial \TT}{\partial t}(x)\right)^2
   + 2  \left(D_y^2 \bfa(\TT x,t)(D\TT(x))_t \right)
     \frac{\partial\TT}{\partial t}(x) \\
& & \quad + 2\left(D_y^2 \bfa_t(\TT x,t)D\TT(x)\right)
    \frac{\partial \TT}{\partial t}(x)
  + 2 D_y \bfa_t(\TT x,t) (D\TT(x))_t \\
& & \quad + \left(D_y^2 \bfa(\TT x,t)D\TT(x)\right)
   \frac{\partial^2 \TT}{\partial t^2}(x)
  + D_y \bfa(\TT x,t) (D\TT(x))_{tt}+ D_y \bfa_{tt}(\TT x, t) (D\TT(x)), 
\end{eqnarray*}
and hence   
\begin{eqnarray*} 
& & D_x \bfb_{tt}(\cdot,t) \to D_x^3 \bfa(\cdot,0)\Sv^2
   + 2 \left(D_x^2 \bfa(\cdot,0)D\Sv\right)\Sv + 2 D_x^2 \bfa_t(\cdot,0)\Sv \\ 
& &  \quad + 2 D_x \bfa_t(\cdot,0) D\Sv + D_x^2\bfa(\cdot,0)\Tv
   + D_x \bfa(\cdot,0)D\Tv + D_x\bfa_{tt}(\cdot,0) 
 \end{eqnarray*} 
as $t \to 0$, uniformly on $\overline{\Omega}$, where $D_x^3 \bfa$ denotes the fourth-order tensor which consists of the third derivatives of the elements of $\bfa$. Thus it follows that  
\begin{eqnarray*}
& & \frac{\partial^2\;}{\partial t^2} \ 
   \mathrm{tr}\left[
  D_x \bfb(\cdot, t) (D\TT)^{-1} \right]
    \mathrm{det}(D\TT) \\ 
& & \quad \rightarrow \ \mathrm{tr} \left[ D_x^3 \bfa(\cdot,0)\Sv^2
   + 2 (D_x^2 \bfa(\cdot,0)D\Sv)\Sv
   + 2 D_x^2 \bfa_t(\cdot,0)\Sv \right] \\
&  & \qquad + \mathrm{tr}\left[
    2 D_x \bfa_t(\cdot,0) D\Sv + D_x^2\bfa(\cdot,0)\Tv
   + D_x \bfa(\cdot,0)D\Tv \right] \\
& & \qquad  + \mathrm{tr} \left[D_x\bfa_{tt}(\cdot,0) + 
  D_x\bfa(\cdot,0)(2(D\Sv)^2 - D\Tv)\right] \\
&  & \qquad + \mathrm{tr} \left[D_x\bfa(\cdot,0)\right] (\nabla\cdot\Tv
   + (\nabla\cdot\Sv)^2 - D\Sv^T : D\Sv) \\
& & \qquad + 2 \mathrm{tr} \left[(D_x^2 \bfa(\cdot,0)\Sv + D_x\bfa(\cdot,0)D\Sv + 
     D_x \bfa_t(x,0))(-D\Sv)\right] \\
& & \qquad + 2 \mathrm{tr} \left[D_x^2 \bfa(\cdot,0)\Sv + D_x\bfa(\cdot,0)D\Sv + 
     D_x \bfa_t(\cdot,0)\right](\nabla\cdot \Sv) \\
&  & \qquad + 2 \mathrm{tr}
    \left[D_x \bfa(\cdot,0) (-D\Sv)\right](\nabla\cdot \Sv)  \\     
& & \quad = \mathrm{tr} \left[ D_x^3 \bfa(\cdot,0)\Sv^2
   + 2 D_x^2 \bfa_t(\cdot,0)\Sv + D_x^2\bfa(\cdot,0)\Tv
    + D_x\bfa_{tt}(\cdot,0)
   \right] \\
& & \qquad + \mathrm{tr} \left[D_x\bfa(\cdot,0)\right] (\nabla\cdot\Tv
   + (\nabla\cdot\Sv)^2 - D\Sv^T : D\Sv) \\
& & \qquad + 2 \mathrm{tr} \left[D_x^2 \bfa(\cdot,0)\Sv 
     + 2D_x \bfa_t(\cdot,0)\right](\nabla\cdot \Sv)  
\end{eqnarray*} 
as $t\rightarrow 0$, uniformly on $\overline{\Omega}$. Hence we obtain  
\begin{eqnarray*} 
& & \left. \frac{\partial^2\;}{\partial t^2} 
   \int_\Omega \mathrm{tr}\left [
  D_x \bfb(\cdot, t) (D\TT)^{-1} \right]
    \mathrm{det}(D\TT) \ \dd x\right\vert_{t=0} \\ 
& & \quad= \int_\Omega \nabla_x\cdot\bfa_{tt}(\cdot,0) + 
   \nabla\cdot [(\nabla\cdot\bfa(\cdot,0))\Tv]
     + 2\nabla\cdot [(\nabla\cdot\bfa(\cdot,0))\Sv] \\
& & \qquad + \mathrm{tr} \left[ D_x^3 \bfa(\cdot,0)\Sv^2
    \right] + 2 
    \mathrm{tr}\left[D_x \bfa_t(\cdot,0)\right](\nabla\cdot \Sv) \\
&   & \qquad + \mathrm{tr} \left[D_x\bfa(\cdot,0)\right] (
    (\nabla\cdot\Sv)^2 - D\Sv^T : D\Sv) \ \dd x,
\end{eqnarray*}
using 
\begin{eqnarray*}
& & \mathrm{tr}\left[D_x^2 \bfa(\cdot,0)\Tv\right]
    + \mathrm{tr}\left[D_x\bfa(\cdot,0)\right](\nabla\cdot\Tv)  
   = \nabla\cdot \left[(\nabla\cdot\bfa(\cdot,0))\Tv
     \right]  \\
& & 2 \mathrm{tr}\left[D_x^2 \bfa_t(\cdot,0)\Sv\right]
    + \mathrm{tr}\left[D_x\bfa_t(\cdot,0)\right](\nabla\cdot\Sv) 
    = 2\nabla\cdot \left[(\nabla\cdot\bfa_t(\cdot,0))\Sv
   \right].
\end{eqnarray*}

We simplify the the last three terms further by the divergence formula. Wright 
\[ \bfa(\cdot,0)=(a^1, \cdots, a^d)^T, \quad \Sv=(S^1, \cdots, S^d)^T, \quad \bnormal=(\nu^1, \cdots, \nu^d)^T, \] 
and 
\[ f_i=\frac{\partial f}{\partial x_i}, \quad f_{ij}=\frac{\partial^2f}{\partial x_i\partial x_j}, \quad f_{ijk}=\frac{\partial^3f}{\partial x_i\partial x_j\partial x_k} \] 
for simplicity. Then it follows that
\begin{eqnarray*}
X : & = & \int_\Omega \mathrm{tr} \left[ D_x^3  \bfa(\cdot,0)\Sv^2
  \right] \dd x  = \sum_{i,p,q}\int_\Omega a^i_{ipq}S^p S^q  \\
& = & - \sum_{i,p,q}\int_\Omega a^{i}_{ip}(S_{q}^p S^q + S^pS_{q}^q)
   + \sum_{i,p,q} \langle a^i_{ip} S^p, 
   S^q \nu^q\rangle_{\partial\Omega},
\end{eqnarray*}
\[ 
 Y := 2 \int_\Omega  \mathrm{tr} \left[D_x^2 \bfa(\cdot,0)\Sv \right] (\nabla\cdot \Sv)
  \ \dd x  = 2 \sum_{i,p,q}\int_\Omega a^i_{ip} S^p S_{q}^q,
\] 
and 
\begin{eqnarray*}
 Z : & = & \int_\Omega \mathrm{tr} \left[D_x\bfa(\cdot,0)\right] (
   (\nabla\cdot\Sv)^2 - D\Sv^T : D\Sv) \ \dd x \\ 
   & = & \sum_{i,p,q}\int_\Omega a^{i}(S_{p}^p S_{q}^q -
    S_{x_p}^qS_{x_q}^p) \ \dd x \\
& = &  \sum_{i,p,q} \int_\Omega -(
   a^i_{ip} S^pS_{q}^q + a^{i}S^pS_{pq}^q) + (a^i_{ip} S^qS_{q}^p + a^{i}_iS^qS_{pq}^p) \ \dd x \\ 
& & + \sum_{i,p,q} \left\{ \langle a^i_i 
    S_{q}^q, S^p \nu^p\rangle_{\partial\Omega} 
- \langle a^i_iS^q,S_{q}^p \nu^p \rangle_{\partial\Omega}\right\} 
\end{eqnarray*}
if $S\in C^{1,1}(\overline{\Omega})$. Adding these equalities, we obtain
\begin{eqnarray*}
 X + Y + Z  & = & \sum_{i,p,q} \left\{ \langle a^i_{ip} S^p,
     S^q \nu^q \rangle_{\partial\Omega}
   + \langle a^{i}S_{q}^q, S^p \nu^p
     \rangle_{\partial\Omega} - \langle
    a^i_iS^q, S_{q}^p\nu^p \rangle_{\partial\Omega} \right\} \\ 
 & = & \left\langle \nabla\cdot
   [(\nabla\cdot\bfa(\cdot,0))\Sv], \Sv\cdot\bnormal 
     \right\rangle_{\partial\Omega} 
   - \left\langle \nabla \cdot \bfa(\cdot,0),
     [(\Sv\cdot \nabla)\Sv]\cdot \bnormal
      \right\rangle_{\partial\Omega}, 
\end{eqnarray*}
which is valid for $S\in C^{0,1}(\overline{\Omega})$. 

Gathering all equations, we obtain the result by the divergence theorem. 
\end{proof}

If $\Omega\subset\R^d$ is a bounded $C^{1,1}$ domain and $\{ \TT\}$ is differentiable, there arises that $\bnormal\in C^{0,1}(\Gamma; S^{d-1})$. We have also $\bfs_i\in C^{0,1}(\Gamma; S^{d-1})$, $1\leq i\leq d-1$, such that 
\[ \{ \bfs_1(\cdot,t), \cdots, \bfs_{d-1}(\cdot,t), \bnormal(\cdot,t)\} \] 
forms a frame of $\partial \Omega_t$ for each $t$. The following lemma ensures again (\ref{by}) for the case that $\Omega$ is $C^{1,1}$. 

\begin{lemma}\label{lemma}
If $\Omega$ is $C^{1,1}$ and $\{ \TT\}$ is differentiable, it holds that 
\[ \bnormal_t=-\sum_{i=1}^{d-1}\left[ \frac{\partial}{\partial \bfs_i}\left(\frac{\partial \TT}{\partial t}\cdot \bnormal\right)\right]\bfs_i, \quad \mbox{a.e. on $\Gamma$}.  \] 
\end{lemma} 

\begin{proof} 
We may fix $x_0\in \Omega$ and assume that $\{ \bfs_1, \cdots, \bfs_{d-1}, \bnormal\}$ are differentiable at $(x,t)=(x_0,0)$ to show the desired equality at this $(x_0,0)$. Write $\bnormal(t)=\bnormal (x_0,t)$, $\bfs_i=\bfs(x_0,0)$, $1\leq i\leq d-1$, and $\bnormal=\bnormal(x_0,0)$. We take the exponential mapping aroud  $x_0$: 
\begin{equation} 
\xi_1\bfs_1 + \cdots + \xi_{d-1}\bfs_{d-1}\in T_{x_0}(\partial \Omega) \ 
  \mapsto  \ x(\bfxi)\in \partial\Omega, \quad \bfxi= (\xi_1, \cdots, \xi_{d-1})\in \R^{d-1}. 
   \label{exponential}
\end{equation} 
This mapping is defined for $\vert \bfxi \vert\ll 1$,  and satisfies $x(0)=x_0$.  Furthermore, it is a local $C^{1,1}$ diffeomorphism, and there arises that 
\[ 
   \frac{\partial x}{\partial\xi_i}
  \biggm|_{\bfxi=\mathbf{0}} = \bfs_i. 
\] 
The perturbed boundary $\partial\Omega_t$ around $x_0(t)=\TT x_0$ is thus 
parametrized by $\bfxi$ as $\TT(x(\bfxi))$, and furthermore, the tangent space $T_{x_0(t)}(\partial\Omega_t)$ is spanned by
\[ 
  \{\tilde\bfs_1(t), \cdots, \tilde\bfs_{d-1}(t)\}, \quad
   \tilde\bfs_i(t)= \left. \frac{\partial\;}{\partial\xi_i}
   \TT(x(\bfxi))\right\vert_{\bfxi=0}, \quad
   i = 1, \cdots, d-1, 
\] 
although $\{ \tilde\bfs_1(t), \cdots, \tilde\bfs_{d-1}(t), \bnormal(t)\}$ does not necessarily form a frame at $x_0(t)=\TT x_0\in \partial \Omega_t$. 

Since $\Omega$ is $C^{1,1}$ and $\{ \TT\}$ is differentiable, these vectors are Lipschitz continuous in $t$, and it follows that 
\[ \bfs_i(t)\cdot\bnormal(t) = 0, \quad \bnormal(t)\cdot\bnormal(t) = 0, \quad \vert t\vert \ll 1. \] 
Then we obtain  
\begin{equation} 
\left.\bnormal_t\right\vert_{t=0}\cdot \bfs_i +
   \bnormal\cdot \left. \frac{\partial \bfs_i}{\partial t}
   \right\vert_{t=0} = 
    \left.\bnormal_t\right\vert_{t=0}
   \cdot\; \bnormal = 0, 
 \label{33x}
\end{equation} 
and furthermore,  
\[ 
  \left. \frac{\partial\bfs_i}{\partial t}\right\vert_{t=0} =
  \left. \frac{\partial^2\;}{\partial t\partial \xi_i}
  \TT(x(\bfxi))\right\vert_{t=0, \ \bfxi = \mathbf{0}} 
  = \left. \frac{\partial x}{\partial \xi_i}S(x(\xi)\right\vert_{\xi=0}=\nabla S\cdot \bfs_i=\frac{\partial}{\partial \bfs_i}\left. \frac{\partial\TT}{\partial t}\right\vert_{t=0}   
\] 
by (\ref{T-taylor}). 

Hence it follows that  
\[ \bnormal_t   = -\sum_{i=1}^{d-1}
\left(\bnormal\cdot \frac{\partial \bfs_i}{\partial t}\right) \bfs_i =-\sum_{i=1}^{d-1}\left[\frac{\partial}{\partial \bfs_i}\left(\frac{\partial \TT}{\partial t}\cdot \bnormal\right)\right]\bfs_i  \] 
at $(x,t)=(x_0,0)$. 
\end{proof} 

\begin{theorem}[second area formula] \label{liouville-area-second}
If $\Omega$ is $C^{2,1}$, $\{\TT\}$ is twice differentiable, $c\in C^{1,1}(\Gamma)$, and $c_{tt}$ is continuous on $\Gamma$, it holds that
\begin{eqnarray*}
& & \left.\frac{\dd^2 \;}{\dd t^2}\int_{\partial \Omega_t}c \ \dd s_t \right\vert_{t=0}
 = \int_{\partial\Omega} \ddot c \ \dd s - \left\langle c_0,  |\nabla_\tau \delta\rho|^2\right\rangle_{\partial\Omega} \\
&  & \quad + \left\langle 
       -2(\Delta_\tau\delta\rho)c_0+2(\nabla\cdot \bnormal)\dot c
     + \nabla\cdot \left[
   ( (\nabla\cdot\bnormal) c_0
      + \frac{\partial c_0}{\partial\bnormal})
    \Sv\right], 
  \delta\rho\right\rangle_{\partial\Omega} -\langle \nabla_\tau^2c_0, (\delta\rho)^2 \rangle_{\partial \Omega} \\
&  & \quad + \left\langle (\nabla\cdot\bnormal) c_0
       + \frac{\partial c_0}{\partial\bnormal},
  \delta^2\rho - ((\Sv\cdot\nabla)\Sv)
   \cdot\bnormal \right\rangle_{\partial\Omega},
\end{eqnarray*}
where $c_0=c(\cdot,0)$, $\dot c=c_t(\cdot, 0)$, $\ddot c=c_{tt}(\cdot,0)$, and 
\[ \Delta_\tau=\nabla_\tau\cdot \nabla_\tau, \quad   
   \nabla_\tau=\left(\frac{\partial}{\partial s_1}, \cdots, \frac{\partial}{\partial s_{d-1}}\right)^T.
\] 
\end{theorem}
\begin{proof}
By the assumption it holds that $\bnormal \in C^{1,1}(\Gamma; S^{d-1})$. Then we apply Lemma~\ref{vector-second} to $\bfa= \bnormal c$.  First, it follows that 
\[ 
\bfa_{tt}= \bnormal_{tt} c+ 2\bnormal_tc_t+ \bnormal c_{tt}, \quad 
\bfa_{t} = \bnormal_{t}c+ \bnormal c_{t}. 
\] 
Second, $\vert \bnormal \vert^2=1$ in $\tilde \Gamma$ implies 
\[ \bnormal\cdot\bnormal_{tt}=-\vert \bnormal_t\vert^2 \quad \mbox{a.e. on $\Gamma$} \] 
as well as (\ref{by}) everywhere. Then we obtain 
\[ \bnormal\cdot \bfa_{tt}= -|\bfal|^2 c+ c_{tt} \] 
almost everywhere with the continuity of its right-hand side on $\Gamma$, where 
\[ \bfal = \left(\frac{\partial\TT}{\partial\bfs_1}
   \cdot\bnormal, \cdots,
   \frac{\partial\TT}{\partial\bfs_{d-1}}
   \cdot\bnormal\right)^T. \] 
It holds also that  
\begin{eqnarray*}
\nabla\cdot\bfa_t & = & [(\nabla\cdot\bnormal) c+\bnormal\cdot \nabla c]_t 
 =   [(\nabla\cdot\bnormal) c]_t + 
   \bnormal_t \cdot \nabla c
   + \bnormal\cdot\nabla c_t\\
& = & [(\nabla\cdot\bnormal) c]_t + 
   \sum_{i=1}^{d-1}\left(\frac{\partial\Sv}{\partial\bfs_i}
   \cdot\bnormal\right) \frac{\partial c}{\partial\bfs_i}
   + \frac{\partial c_t}{\partial\bnormal} \\
& = & [
     (\nabla\cdot\bnormal) c]_t+
     \bfal\cdot\nabla_\tau c
  + \frac{\partial c_t}{\partial\bnormal} 
\end{eqnarray*}
almost everywhere with its right-hand side on $\Gamma$. Hence Lemma \ref{vector-second} is applicable. 

Since $\Omega$ is $C^{2,1}$ and $\{ \TT\}$ is twice differentiable, $\bnormal_t$ in Lemma \ref{lemma} is Lipschitz continuous on $\Gamma$, and it holds that 
\[ \nabla\cdot \bnormal_t=-\sum_{i=1}^{d-1}\frac{\partial^2}{\partial \bfs_i^2}\left(\frac{\partial\TT}{\partial t}\cdot \bnormal\right).  \] 
We thus obtain 
\begin{eqnarray*} 
\left. [(\nabla\cdot \bnormal)c\right]_t & = & (\nabla\cdot \bnormal_t)c+(\nabla\cdot \bnormal)c_t \\ 
& = & - \left(\Delta_\tau \frac{\partial \TT}{\partial t}\right)c+(\nabla\cdot \bnormal)c_t. 
\end{eqnarray*} 
 
Then we obtain the result by 
\[ \nabla\cdot\bfa = (\nabla\cdot\bnormal) c
       + \frac{\partial c}{\partial\bnormal}, \quad 
       \left. \left\vert \alpha\right\vert^2\right\vert_{t=0}=\sum_{i=1}^{d-1}\left(\frac{\partial \delta\rho}{\partial s_i}\right)^2=\vert \nabla_\tau \delta \rho\vert^2, 
\]  
and 
\[ 
2\left. \left\langle
\alpha \cdot \nabla_\tau c, \frac{\partial \TT\cdot \bnormal}{\partial t} \right\rangle_{\partial\Omega}\right\vert_{t=0}
=\langle \nabla_\tau c_0, \nabla_\tau (\delta \rho)^2\rangle= -\langle \nabla_\tau^2c_0, (\delta \rho)^2 \rangle    
\] 
derived from 
\[ \bfal\cdot \nabla_\tau c=\sum_{i=1}^{d-1}\frac{\partial\TT\cdot \bnormal}{\partial s_i} \frac{\partial c}{\partial \bfs_i} =\nabla_\tau\frac{\partial \TT\cdot \bnormal}{\partial t}\cdot \nabla_\tau c. \] 
\end{proof}

\section{Hadamard variation of the Green's function}\label{sec4}

\subsection{First variation} 

The existence of Hadamard variation (\ref{6}) on the Green's function is assured by the method of \cite{SuzTsu16}. Recall that $N(x,y,t)$ and $N(x,y)$ denote the Green's function on $\Omega_t$ and $\Omega$ defined by (\ref{4.0}), (\ref{4.1}), (\ref{4.5}) and (\ref{4.0}), (\ref{4}), (\ref{5}), using $u_t\in H^1(\Omega_t)$ and $u\in H^1(\Omega)$, respectively. 

\begin{theorem}\label{thm15} 
Let $\Omega\subset \R^d$ be a bounded Lipschitz domain and fix $y\in \Omega$. Let $u=u(\cdot,t)\in H^1(\Omega_t)$ be the solution to (\ref{4.5}) in Theorem \ref{thm0}. Then, if $\{ \TT\}$ is differentiable, it holds that  
\begin{equation} 
v\in C^1(-\varepsilon, \varepsilon; H^1(\Omega)), \quad v(x,t)=u(\TT x, t). 
 \label{44}
\end{equation} 
In particular, the first Hadamard variation in (\ref{6}), 
\[ \delta N(\cdot,y)=\left. \frac{\partial u}{\partial t}(\cdot,t)\right\vert_{t=0}, \] 
exists in the sense of distributions in $\Omega$. There arises that 
\begin{equation} 
\Delta \delta N(\cdot,y)=\Delta \dot u=0 \quad \mbox{in $\Omega$}, 
 \label{47.1}
\end{equation} 
and furthermore, $\delta N(\cdot,y)\in L^2(\Omega)$, more precisely,  
\begin{equation} 
\delta N(x,y)=\left. \frac{\partial v}{\partial t}(x,t)\right\vert_{t=0}-(S\cdot\nabla)u_0(x) 
 \label{45}
\end{equation}
for $S\in C^{0,1}(\overline{\Omega}, \R^d)$ defined by (\ref{taylor0}) and $u_0=u(\cdot,0)\in H^1(\Omega)$.  
\end{theorem} 

\begin{proof} 
To show (\ref{44}), we take $\psi\in C_0^\infty(\R^d)$, $0\leq \psi \leq 1$, such that $\psi=0$ and $\psi=1$ in $B(y,r)$ and $\R^d\setminus B(y,2r)$, respectively, for $0<r\ll 1$. Then $\tilde \Gamma=\Gamma(\cdot-y)\varphi$ is independent of $t$, and $w=u+\tilde \Gamma$ satisfies  
\begin{equation} 
-\Delta w=h \ \mbox{in $\Omega_t$}, \quad w=0 \ \mbox{on $\gamma_t^0$}, \quad \frac{\partial w}{\partial \bnormal}=0 \ \mbox{on $\gamma_t^1$}  
 \label{45x}
\end{equation} 
for $h=-\Delta \tilde \Gamma$. The Poisson problem (\ref{45x}) takes the weak form 
\[  
w\in V_t, \quad  
\int_{\Omega_t}\nabla w\cdot \nabla \varphi \ \dd x=\int_{\Omega_t}h\varphi \ \dd x, \ \forall \varphi\in V_t,   
\]  
where $V_t=\{ v\in H^1(\Omega_t) \mid \left. v\right\vert_{\gamma_t^1}=0\}$. 
Then we obtain 
\[  z\in C^1(-\varepsilon, \varepsilon; H^1(\Omega)), \quad z(x,t)=w(\TT x, t) \] 
and hence (\ref{44}) as in \cite[Theorem 16]{SuzTsu16}, because $\{ \TT\}$ is differentiable.  

Having (\ref{44}), the rest part of this theorem follows as in \cite{SuzTsu16}
\end{proof} 

If $\Omega$ is $C^{1,1}$, furthermore, we have $u_0\in H^2(\Omega)$ by the elliptic regularity, recalling (\ref{assum-gamma}).  Then it holds that $\delta N(\cdot,y)\in H^1(\Omega)$ by (\ref{44})-(\ref{45}), and then $H^1$ theory is applicable to (\ref{4}) for $u=\delta N(\cdot,y)$. We have also $N(\cdot, y)\in H^2(\Omega)$, which implies  
\begin{equation} 
\left.\nabla N(\cdot,y)\right\vert_{\partial \Omega}\in H^{1/2}(\partial \Omega).  
 \label{41}
\end{equation} 
Thus we have the well-definedness of the right-hand side of the desired identity in the following theorem.  

\begin{theorem}[first variational formula]
\label{firstHadamard}
If $\Omega$ is $C^{1,1}$ and $\{\TT\}$ is differentiable, it holds that 
\[ 
\delta N(x,y) = 
\left\langle \delta\rho\frac{\partial N}{\partial \bnormal}(\cdot,x),
\frac{\partial N}{\partial \bnormal}(\cdot,y)
\right\rangle_{\gamma^0} 
 - \left\langle \delta\rho \nabla_{\tau} N(\cdot,x), \nabla_{\tau} N(\cdot,y) \right\rangle_{\gamma^1}, \quad x,y \in \Omega, 
\] 
where $\delta\rho= \Sv\cdot\bnormal$ and $\nabla_\tau$ is the
tangential gradient on $\partial\Omega$ defined by \eqref{tan-grad}.
\end{theorem}

\begin{proof}
Continue to fix $y\in \Omega$. We have readily confirmed $N(\cdot,y)\in H^2(\Omega)$ and $\dot u=\delta N(\cdot,y)\in H^1(\Omega)$. Now we show that $\dot u=\delta N(\cdot,y)\in H^1(\Omega)$ solves 
\begin{equation} 
\Delta \dot u=0 \ \mbox{in $\Omega$}, \quad \dot u =-\delta\rho\frac{\partial N}{\partial\bnormal}(\cdot, y) \ \mbox{on $\gamma^0$}, \quad \frac{\partial \dot u}{\partial \bnormal}=\nabla_\tau\cdot(\delta\rho\nabla_\tau N(\cdot,y)) \ \mbox{on $\gamma^1$}, 
 \label{49}
\end{equation} 
where 
\[ \nabla_\tau\cdot(\delta \rho\nabla_\tau N(\cdot,y))=\sum_{i=1}^{d-1}\frac{\partial}{\partial \bfs_i}\left(\delta \rho \frac{\partial N}{\partial \bfs_i}N(\cdot, y)\right). \]  

Once (\ref{49}) is shown, Theorem \ref{thm0} is applicable to this Poisson equation, because (\ref{41}) implies 
\begin{equation} 
\delta\rho\frac{\partial N}{\partial\bnormal}(\cdot, y)\in H^{1/2}(\gamma^0), \quad  
\nabla_\tau\cdot(\delta \rho\nabla_\tau N(\cdot,y)) \in H^{-1/2}(\gamma^1)  
 \label{43}
\end{equation} 
by $\delta \rho\in C^{0,1}(\partial \Omega)$. Then the desired equality follows from the representation formula (\ref{8}) of the solution $\dot u=\dot u(x)$ to (\ref{49}), because $N(x,y)=N(y,x)$.  
 
Since the boundary condition of $v=\delta N(\cdot,y)$ on $\gamma^0$ in (\ref{49}) is assured by the result in \cite{SuzTsu16} on the Dirichlet boundary condition, we have only to confirm the boundary condition on $\gamma^1$ in (\ref{43}). To this end, we take an open neighbourhood of $\gamma^1$, denoted by $\tilde \Omega$, satisfying $\tilde \Omega\cap \gamma^0=\emptyset$ and $y\not\in \tilde \Omega$. 

Let $\varphi\in C_0^\infty(\tilde \Omega)$. Then, for $\vert t\vert\ll 1$ it holds that 
\begin{equation}
  \int_{\Omega_t}
  \nabla_x N(x,y,t)\cdot\nabla\varphi(x)  \ \dd x =0  
  \label{green-one}
\end{equation}
by 
\[ \Delta N(\cdot,y,t)=0 \ \mbox{in $\Omega_t\setminus \{y\}$}, \quad \frac{\partial N}{\partial \bnormal}(\cdot,y,t)=0 \ \mbox{on $\gamma_t^1$},  \] 
where the normal derivative of $N(\cdot,y,t)$ on $\gamma_t^1$ belongs to $H^{-1/2}(\gamma_t^1)$. 

We apply Liouville's first volume formula, Theorem \ref{thm8}, to \eqref{green-one} for 
\begin{equation} 
c(x,t) = \nabla_x N(x,y,t)\cdot\nabla  \varphi(x),  
 \label{c1}
\end{equation} 
satisfying $c\in L^\infty(Q)$ and $\nabla c\in L^\infty(Q;\R^d)$.  In fact, we obtain  
\[ b\in C^1(-\varepsilon, \varepsilon; L^2(\Omega)), \quad \beta \in C^0(-\varepsilon, \varepsilon; L^2(\Omega)) \] 
for $b(x,t)=c(\TT x, t)$ and $\beta(x,t)=\nabla_x c(\TT x, t)\cdot \frac{\partial \TT}{\partial t}(x)$ by Theorem \ref{thm15}.  Then it follows that 
\[ \alpha\in C^0(-\varepsilon, \varepsilon; L^2(\Omega)) \] 
for $\alpha(x,t)=c_t(\TT x,x)$ from (\ref{35x}). 

Hence there holds that   
\begin{eqnarray*}
  0 & = & \left. \frac{\dd\;}{\dd t}
    \int_{\Omega_t} \nabla_x N(x,y,t)\cdot\nabla\varphi(x) \ 
 \dd x \right\vert_{t=0} \\ 
 & = & \int_\Omega \dot c+\nabla\cdot (c_0S) \ dx \\ 
    & = & \int_{\Omega}
    \nabla\dot u(x) \cdot\nabla \varphi(x) 
   + \nabla_x\cdot \left(\nabla_xN(x,y)\cdot\nabla \varphi(x)\frac{\partial \TT}{\partial t}(x)\right) \ \dd x \\ 
   & = & \left\langle \varphi, \frac{\partial\dot u}{\partial \bnormal}\right\rangle_{\gamma^1} + \int_{\gamma^1}\delta \rho \nabla_xN(\cdot, y)\cdot \nabla \varphi \ ds \\
   & = & \left\langle \varphi, \frac{\partial\dot u}{\partial \bnormal}\right\rangle_{\gamma^1}+ \int_{\gamma^1}\delta \rho \nabla_\tau N(\cdot, y)\cdot \nabla_\tau\varphi \ ds \\ 
  & = & \left\langle \varphi, \frac{\partial \dot u}{\partial \bnormal}\right\rangle_{\gamma^1}-\langle \varphi, \, \nabla_\tau(\delta \rho\nabla_\tau N(\cdot,y)) \rangle_{\gamma^1} 
\end{eqnarray*}
by $\varphi \in C_0^\infty(\tilde \Omega)$ and (\ref{47.1}). Then we obtain 
\[ \frac{\partial\dot u}{\partial \bnormal}=\nabla_\tau\cdot (\delta\rho \nabla_\tau N(\cdot,y)) \quad \mbox{on $\gamma^1$} \] 
as an element in $H^{-1/2}(\gamma^1)$, because $\varphi\in C_0^\infty(\tilde \Omega)$ is arbitrary. 
\end{proof}

\subsection{Second variation} 

We begin with the existence of the second variation $\delta^2N$ in (\ref{77}) as in Theorem \ref{thm15}. 

\begin{theorem}
Let $\Omega\subset \R^d$ be a bounded Lipschitz domain, $\{ \TT\}$ be twice differentiable, and $u\in H^1(\Omega_t)$ be the solution to (\ref{4.5}) for $y\in \Omega$. Then it holds that 
\begin{equation} 
v\in C^2(-\varepsilon, \varepsilon; H^1(\Omega)),  \quad v(x,t)=u(\TT x,t).  
 \label{54}
\end{equation} 
In particular, $\ddot u=\delta^2N(\cdot,y)$ in (\ref{77}) exists in the sense of distributions in $\Omega$, and there arises that 
\[ \Delta \ddot u=0 \quad \mbox{in $\Omega$}. \] 
If $\Omega$ is $C^{1,1}$ and $C^{2,1}$, furthermore, this $\ddot u$ belongs to $L^2(\Omega)$ and $H^1(\Omega)$, respectively. 
\end{theorem} 

\begin{proof} 
All the results except for the regularity of $\ddot u$ follow from the weak form (\ref{45x}) as in \cite{SuzTsu16}. There arises also that 
\begin{equation} 
\left. \frac{\partial^2v}{\partial t^2}\right\vert_{t=0}=\ddot u+2S\cdot \nabla \dot u+R\cdot \nabla u+[\nabla^2u]S\cdot S\in H^1(\Omega).
 \label{55}
\end{equation} 

If $\Omega$ is $C^{1,1}$ we have $u\in H^2(\Omega)$ and hence $\dot u\in H^1(\Omega)$ by (\ref{45}). Then $\ddot u\in L^2(\Omega)$ follows from (\ref{55}) and $S, R\in C^{0,1}(\overline{\Omega})$. If $\Omega$ is $C^{2,1}$ there arises $u\in H^3(\Omega)$ and hence $\dot u\in H^2(\Omega)$ by (\ref{49})  and (\ref{5}). Then (\ref{55}) implies $\ddot u\in H^1(\Omega)$, similarly. 
\end{proof} 

\begin{lemma}\label{lem19} 
If $\Omega$ is $C^{2,1}$ and $\{ \TT\}$ is twice differentiable, $\ddot u=\delta^2N(\cdot,y)\in H^1(\Omega)$ satisfies 
\begin{equation} 
\Delta \ddot u=0 \ \mbox{in $\Omega$}, \quad \ddot u= g \ \mbox{on $\gamma^0$}, \quad 
\frac{\partial \ddot u}{\partial\bnormal}= h \quad \mbox{on $\gamma^1$} 
 \label{58}
\end{equation} 
for 
\begin{eqnarray} 
& & g= - \chi
\frac{\partial N}{\partial\bnormal}(\cdot,y)
+ 2 \delta\rho\frac{\partial \dot u}{\partial\bnormal} \nonumber\\ 
& & h= \nabla_\tau\cdot (\sigma N_\tau(\cdot,y))+2\nabla_\tau\cdot (\delta\rho\nabla_\tau \dot u), 
 \label{58.1}
\end{eqnarray} 
where $\delta \rho=S\cdot\bnormal$ and  
\begin{eqnarray} 
& & \chi = \delta^2\rho+\delta\rho\frac{\partial\delta\rho}{\partial\bnormal}+(\nabla\bnormal)[S,S]-(\delta\rho)^2\nabla\cdot\bnormal -(S\cdot\nabla)\delta \rho \nonumber\\
& & \sigma = \delta^2\rho-2(S_\tau\cdot \nabla_\tau)\delta\rho+(\nabla \bnormal)[S,S],   
  \label{59} 
\end{eqnarray} 
where $S_\tau=S-(\delta\rho)\bnormal$. 
\end{lemma}

\begin{proof}
We have readily obtained  
\begin{equation} 
u=N(\cdot,y)-\Gamma(\cdot-y)\in H^3(\Omega), \ \dot u=\delta N(\cdot,y)\in H^2(\Omega), \ \ddot u=\delta^2N(\cdot,y)\in H^1(\Omega)  
 \label{62}
\end{equation} 
if $\Omega$ is $C^{2,1}$ and $\{ \TT\}$ is twice differentiable. By the same reason there arises that 
\[ \chi, \sigma \in C^{0,1}(\partial \Omega) \] 
for $\chi$ and $\sigma$ in (\ref{59}). Hence it follows that 
\[ g\in H^{1/2}(\gamma^0), \quad h\in H^{-1/2}(\gamma^1) \] 
for $g$ and $h$ defined by (\ref{58.1}). 

We have confirmed $\Delta \ddot u=0$ in $\Omega$ in the previous theorem. It is also shown that 
\[ \ddot u=g=-\chi\frac{\partial N}{\partial\bnormal}(\cdot,y)+2\delta\rho\frac{\partial\dot u}{\partial \bnormal} \quad \mbox{on $\gamma^0$} \] 
for 
\begin{equation} 
\chi=(R-(S\cdot \nabla)S)\cdot \bnormal-(\delta\rho)^2(\nabla\cdot \bnormal)-(S\cdot \nabla)\delta\rho+\frac{\partial(\delta\rho)^2}{\partial \bnormal} 
 \label{xx3}
\end{equation} 
by \cite{SuzTsu16}. Then there arises the first equality of (\ref{59}) by $R\cdot \bnormal=\delta^2\rho$,  
\begin{eqnarray} 
[(S\cdot\nabla)S]\cdot \bnormal & = & (S\cdot\nabla)(S\cdot \bnormal)-[(S\cdot \nabla)\bnormal]\cdot S \nonumber\\ 
& = & (S_\tau\cdot \nabla_\tau)\delta\rho+\delta \rho\frac{\partial \delta\rho}{\partial \nu}-S\cdot [(S\cdot \nabla)\bnormal], 
 \label{xx1}  
\end{eqnarray} 
and 
\begin{equation} 
S\cdot[(S\cdot\nabla)\bnormal]=(\nabla \bnormal)[S,S]. 
 \label{xx2}
\end{equation} 
It thus suffices to ensure 
\begin{equation} 
\frac{\partial \ddot u}{\partial \bnormal}=h \ \mbox{on $\gamma^1$}. 
 \label{60}
\end{equation} 

For this purpose we use the open neighbourhood $\tilde \Omega$ of $\gamma^1$ in the proof of Theorem \ref{firstHadamard} satisfying $\tilde \Omega \cap \gamma^0=\emptyset$ and $y\not\in \tilde \Omega$.  Taking $\varphi\in C_0^\infty(\tilde\Omega)$, we have $b\in C^2(-\varepsilon, \varepsilon; L^2(\Omega))$ for $b=b(x,t)$ defined by (\ref{c1}) and (\ref{32.0}), from the proof of \cite[Theorem 16]{SuzTsu16}. Hence Theorem \ref{lem12} is applicable. 

Since $\Omega$ is $C^{1,1}$ it holds that $\nabla^2c\in L^\infty(Q; \R^d\times \R^d)$. Then we obtain     
\[ b_{tt}=\gamma+2\delta+\mu+\sigma \] 
in the sense of distributions in $\tilde Q=\Omega\times (-\varepsilon, \varepsilon)$, 
where 
\begin{eqnarray*}
& & \gamma(x,t)=c_{tt}(T_tx, t), \qquad \qquad \quad \ \delta(x,t)=\beta_t(x,t)-\nabla c(T_tx,t)\cdot \frac{\partial^2\TT}{\partial t^2}(x) \\ 
& & \sigma(x,t)=\nabla c(T_tx, t)\cdot \frac{\partial^2\TT}{\partial t^2}(x), \quad  
\mu(x,t)=\nabla^2c(T_tx,t)\left[\frac{\partial \TT}{\partial t}(x), \frac{\partial \TT}{\partial t}(x)\right]. 
\end{eqnarray*} 
By (\ref{54}) we have $\gamma, \delta, \sigma \in C(-\varepsilon, \varepsilon; L^2(\Omega))$, which implies $\mu \in C(-\varepsilon, \varepsilon; L^2(\Omega))$. 

Put $t=0$ in the conclusion of Lemma \ref{lem12} and apply (\ref{62}). Then we obtain 
\begin{eqnarray*}
  0 & = & \left. \frac{\dd^2\;}{\dd t^2} 
   \int_{\Omega_t} \nabla N(\cdot,y,t)\cdot\nabla\varphi \ \dd x\right|_{t=0} \\ 
   & = & \int_\Omega \nabla \delta^2 N(\cdot,y)\cdot\nabla \varphi \ \dd x
    + 2 \langle \nabla \varphi, \, \delta\rho
   \nabla\delta N(\cdot,y) \rangle_{\partial\Omega} + \left\langle \nabla\cdot\left[
    (\nabla_x N(\cdot,y)\cdot\nabla \varphi) \Sv
     \right], \, \delta\rho \right\rangle_{\partial\Omega} \\ 
     & & + \left\langle \nabla \varphi, \,
    [ (\Tv - (\Sv\cdot\nabla)\Sv)\cdot\bnormal] \nabla N(\cdot,y)\right\rangle_{\partial\Omega} 
\end{eqnarray*}
by the proof of Theorem \ref{second-liouville-thm}. 

We examine each term on the right-hand side, recalling $\varphi\in C_0^\infty(\tilde \Omega)$.  First, it follows that 
\[ \int_\Omega \nabla\delta^2 N(\cdot,y)\cdot\nabla \varphi \ \dd x
   = \left\langle \varphi,\, \frac{\partial\;}{\partial\bnormal}
      \delta^2 N(\cdot,y) \right\rangle_{\gamma^1}
\] 
from $\Delta \ddot u=0$. Second, we have 
\begin{eqnarray*}
\langle \nabla \varphi, \, \nabla N(\cdot,y) 
     (\Tv - (\Sv\cdot \nabla) \Sv)\cdot\bnormal \rangle_{\partial\Omega} 
& = & \langle
     \nabla_\tau \varphi, 
    [(\Tv - (\Sv\cdot\nabla)\Sv)\cdot\bnormal]
    \nabla_\tau N(\cdot,y) \rangle_{\gamma^1} \\ 
& = & - \langle \varphi, \,
     \nabla_\tau \cdot 
    ([(\Tv - (\Sv\cdot\nabla)\Sv)\cdot\bnormal]
   \nabla_\tau N(\cdot,y)) \,
     \rangle_{\gamma^1}.
\end{eqnarray*}
by \eqref{dif-form} for $F = \varphi$, $g = (\Tv - (\Sv\cdot\nabla)\Sv)\cdot\bnormal$, and 
$H = N(\cdot,y)$. 

Third, there arises that 
\begin{eqnarray*}
& & \left\langle \nabla \cdot [
    (\nabla N(\cdot,y)  \cdot\nabla \varphi) \Sv], \, \delta\rho \right\rangle_{\partial\Omega} \\
& & \quad = \left\langle  \left[
     \sum_{i=1}^{d-1} \frac{\partial\;}{\partial \bfs_i}
     ([\Sv\cdot\bfs_i] \nabla N(\cdot,y)\cdot\nabla \varphi)
     + \frac{\partial\;}{\partial \bnormal}
     ([\Sv\cdot\bnormal] \nabla N(\cdot,y)\cdot\nabla \varphi)\right], \, \Sv\cdot\bnormal\right\rangle_{\gamma^1} \\
& & \quad = \left\langle \nabla N(\cdot,y)\cdot\nabla \varphi, \,
     (\Sv\cdot\bnormal)
     \frac{\partial(\Sv\cdot\bnormal)}{\partial\bnormal}
    - \sum_{i=1}^{d-1}(\Sv\cdot\bfs_i)
     \frac{\partial(\Sv\cdot\bnormal)}{\partial\bfs_i}
     \right\rangle_{\gamma^1} \\
& & \qquad + \left\langle
    \frac{\partial\;}{\partial\bnormal}
    (\nabla N(\cdot,y)\cdot\nabla \varphi), \,
    (\Sv\cdot\bnormal)^2 \right\rangle_{\gamma^1} \\
& & \quad = \left\langle 
    \nabla_{\tau}N(\cdot,y) \cdot \nabla_{\tau}\varphi, \,
     \delta \rho\frac{\partial \delta\rho}{\partial\bnormal}
     - (\Sv\cdot\nabla_\tau)\delta\rho
      \right\rangle_{\gamma^1} 
      + \sum_{i=1}^{d-1} \left\langle 
      \frac{\partial^2 N}{\partial\bnormal^2}
      \frac{\partial\varphi}{\partial\bnormal}
      + \frac{\partial N}{\partial \bfs_i}
       \frac{\partial^2\varphi}{\partial \bfs_i\partial\bnormal}, (\delta \rho)^2
   \right\rangle_{\gamma^1} 
\end{eqnarray*} 
by \eqref{dif-form0} for $F = N(\cdot,y)\cdot\nabla \varphi$ and $H = \Sv\cdot\bnormal$. Here we have 
\begin{eqnarray*} 
& &  \left\langle 
    \nabla_{\tau}N(\cdot,y) \cdot \nabla_{\tau}\varphi, \,
     \delta \rho\frac{\partial\delta \rho}{\partial\bnormal}
     - (\Sv\cdot\nabla_\tau)\delta \rho
      \right\rangle_{\gamma^1} \\ 
& & \quad = - \left\langle \varphi, \,
     \nabla_\tau\cdot [
    (\delta \rho\frac{\partial\delta\rho}{\partial\bnormal}
    - (\Sv\cdot\nabla_{\tau})\delta\rho)
      \nabla_\tau N(\cdot,y)]\right\rangle_{\gamma^1} 
\end{eqnarray*} 
and also 
\begin{eqnarray*} 
& & \left\langle \sum_{i=1}^{d-1} \left(
      \frac{\partial^2 N}{\partial\bnormal^2}
      \frac{\partial\varphi}{\partial\bnormal}+\frac{\partial N}{\partial \bfs_i}
       \frac{\partial^2\varphi}{\partial \bfs_i\partial\bnormal}\right), \,  (\delta\rho)^2
   \right\rangle_{\gamma^1} \\ 
& & = - \left\langle
    \frac{\partial\varphi}{\partial\bnormal},\,
      (\delta\rho)^2(\nabla_\tau)^2 N(\cdot,y)
     \right\rangle_{\gamma^1}
  - \left\langle \frac{\partial \varphi}{\partial \nu}, \nabla_\tau (\delta\rho)^2\cdot\nabla_\tau N(\cdot,y) \right\rangle  \\ 
& &  = - 2 \left\langle \frac{\partial \varphi}{\partial \bnormal}, 
    (\delta\rho)^2(\nabla_\tau)^2N(\cdot,y)\right\rangle_{\gamma^1} 
    - \left\langle
      \frac{\partial\varphi}{\partial\bnormal}, \,
      \nabla_\tau(\delta\rho)^2\cdot \nabla_\tau N(\cdot,y)\right\rangle_{\gamma^1} 
\end{eqnarray*}
by 
\[ 
   \frac{\partial N}{\partial\bnormal}(\cdot,y)
  = \frac{\partial^2 N}{\partial \bfs_i\partial\bnormal}(\cdot,y)
  = 0, \ 
  \sum_{i=1}^{d-1} 
   \frac{\partial^2 N}{\partial \bfs_i^2}(\cdot,y)
   + \frac{\partial^2 N}{\partial\bnormal^2}(\cdot,y) = 0 \quad \mbox{on $\gamma^1$}.  
\] 
Finally, we notice (\ref{49}) to deduce 
\begin{eqnarray*}
& & 2 \langle \nabla \varphi, \, \delta\rho
   \nabla \delta N(\cdot,y)\rangle_{\gamma^1} \\  
& & \quad = 2 \left\langle \frac{\partial \varphi}{\partial \bnormal}, \,
    \delta\rho\frac{\partial \delta N}{\partial \bnormal}(\cdot,y)
     \right\rangle_{\gamma^1} + 2 
     \left\langle \nabla_\tau \varphi, 
      \delta\rho\nabla_\tau \delta N(\cdot,y)  \right\rangle_{\gamma^1} \\ 
&  & \quad =  2 \left\langle
   \frac{\partial \varphi}{\partial \bnormal}, \,
    \delta\rho
   \nabla_\tau(\delta\rho)\cdot \nabla_\tau N(\cdot,y) \right\rangle_{\gamma^1}
   + 2 \left\langle
     \frac{\partial \varphi}{\partial \bnormal}, \,
    (\delta\rho)^2\nabla_\tau^2N(\cdot,y) \right\rangle_{\gamma^1} \\ 
    & & \qquad -  2 \left\langle \varphi, \,
      \nabla_\tau\cdot(\delta\rho\nabla_\tau \delta N(\cdot,y))\right\rangle_{\gamma^1}.
\end{eqnarray*}

Gathering these equalities, we obtain
\begin{eqnarray*}
& & 0 =  \left\langle \varphi, \,
    \frac{\partial\;}{\partial\bnormal}
    \delta^2 N(\cdot,y) \right\rangle_{\gamma^1} - \left\langle  \varphi, \,
     \nabla_\tau\cdot \left[
     \left( \delta\rho
     \frac{\partial \delta \rho}{\partial \bnormal}
      - (\Sv\cdot\nabla_{\tau})\delta\rho \right)
     \nabla_\tau N(\cdot,y)
      \right] \right\rangle_{\gamma^1} \\
& & \quad - \langle \varphi, \nabla_\tau\cdot 
     ([(\Tv - (\Sv\cdot\nabla)\Sv)\cdot\bnormal] \nabla_\tau N(\cdot,y)) \rangle_{\gamma^1} - 2 \langle \varphi, \nabla_\tau \cdot (
   \delta\rho \nabla_\tau \delta N(\cdot,y))\rangle_{\gamma^1} 
\end{eqnarray*}
and hence the result because $\varphi \in C_0^\infty(\tilde \Omega)$ is arbitrary. In fact, we obtain (\ref{58.1}) for 
\begin{equation} 
\sigma = \delta\rho \frac{\partial\delta\rho}{\partial \bnormal}-(S_\tau\cdot\nabla_\tau)\delta\rho +(R-(S\cdot\nabla)S)\cdot \bnormal 
 \label{xx4}
\end{equation} 
and then, the second equality of (\ref{59}) follows from $R\cdot \bnormal=\delta^2\rho$, (\ref{xx1}), and (\ref{xx2}).  
\end{proof} 

\begin{theorem}[second variational formula]\label{thm21}
If $\Omega$ is $C^{2,1}$ and $\{ \TT\}$ is twice differentiable, it holds that 
\begin{eqnarray*}  
& & \delta^2 N(x,y) =  -2 (\nabla \delta N(\cdot,x),
  \nabla \delta N(\cdot,y) )
  + \left\langle \chi \frac{\partial N}{\partial \bnormal}(\cdot,x),
  \frac{\partial N}{\partial \bnormal}(\cdot,y)\right\rangle_{\gamma^0} \\ 
& & \quad - \left\langle \sigma \nabla_{\tau} N(\cdot,x),
   \nabla_{\tau} N(\cdot,y) \right\rangle_{\gamma^1} 
\end{eqnarray*}  
for $x,y \in \Omega$, where $( \ , \ )$ denotes the inner product in $L^2(\Omega)$. 
\end{theorem}
\begin{proof} 
Form Lemma \ref{lem19} and the representation formula (\ref{8}), it follows that 
\begin{equation} 
\delta^2N(x,y)=\ddot u(x)=- \left\langle g, \frac{\partial N}{\partial \bnormal}(\cdot,x)\right\rangle_{\gamma^0}+\langle N(\cdot,x), h \rangle_{\gamma^1} 
 \label{64}
\end{equation} 
for $g$, $h$ defined by (\ref{58.1})-(\ref{59}), where $x,y\in \Omega$. 

By (\ref{49}), we obtain 
\begin{eqnarray*}
   0 & = & \left\langle \frac{\partial \delta N}{\partial\bnormal}(\cdot,x), \,
    \delta N(\cdot,y) + \delta \rho \frac{\partial N}{\partial\bnormal}(\cdot,y)
   \right\rangle_{\gamma^0} \\ 
   & & + \left\langle \delta N(\cdot,y), \, 
   \frac{\partial \delta N}{\partial\bnormal}(\cdot,x)
    - \nabla_\tau\cdot (\delta \rho \nabla_\tau N(\cdot,x))\right\rangle_{\gamma^1} \\ 
   & = & \left\langle  \delta N(\cdot,y), \,
   \frac{\partial \delta N}{\partial\bnormal}(\cdot,x)
     \right\rangle_{\partial\Omega} + \left\langle \frac{\partial \delta N}{\partial\bnormal}(\cdot,x), \delta\rho \frac{\partial N}{\partial\bnormal}(\cdot,y) \right\rangle_{\gamma^0} \\ 
& & + 
      \left\langle \nabla_\tau \delta N(\cdot,y), 
       \delta \rho \nabla_\tau N(\cdot,x) \right\rangle_{\gamma^1}  \\ 
   & = & (\nabla\delta N(\cdot,y) , \nabla\delta N(\cdot,x)) + \left\langle \frac{\partial \delta N}{\partial\bnormal}(\cdot,x), \delta\rho \frac{\partial N}{\partial\bnormal}(\cdot,y)\right\rangle_{\gamma^0}
     \\ 
     & & - \langle N(\cdot,x), 
       \nabla_\tau\cdot (\delta \rho\nabla_\tau \delta N(\cdot,y)) \rangle_{\gamma^1}
\end{eqnarray*}
for $x\in \Omega$, and hence 
\begin{eqnarray}
   (\nabla\delta N(\cdot,x), \nabla\delta N(\cdot,y)) & = & - \left\langle \delta\rho \frac{\partial N}{\partial\bnormal}(\cdot,y), \frac{\partial \delta N}{\partial\bnormal}(\cdot,x) \right\rangle_{\gamma^0} \nonumber\\ 
   & & + \langle N(\cdot,x), \nabla_\tau\cdot (\delta \rho \nabla_\tau \delta N(\cdot,y)) \rangle_{\gamma^1}.
\label{lem-delta-N}
\end{eqnarray}
Then the result follows from (\ref{64})-(\ref{lem-delta-N}) as 
\begin{eqnarray*} 
& & \delta^2N(x,y) - \left\langle\chi \frac{\partial N}{\partial \nu}(\cdot,x), \frac{\partial N}{\partial \nu}(\cdot,y) \right\rangle_{\gamma^0}+\langle \sigma \nabla_\tau N(\cdot,x), \nabla_\tau N(\cdot,y) \rangle_{\gamma^1} \\ 
& & \quad = -2 \left\langle \delta \rho \frac{\partial N}{\partial \nu}(\cdot,x), \frac{\partial N}{\partial \nu}(\cdot,y) \right\rangle_{\gamma^0}+2\langle N(\cdot,x), \nabla_\tau(\delta \rho \nabla_\tau\delta N(\cdot,y) \rangle_{\gamma^1} \\ 
& & \quad = -2(\nabla \delta N(\cdot, x), \nabla \delta N(\cdot,y) ).  
\end{eqnarray*}

\end{proof} 

\appendix 

\section{Liouville's area formula and differential forms}\label{seca2}

There is an alternative argument for the proof of Lemmas \ref{vector-first} and \ref{vector-second} using differential forms. Here we assume that $\Omega\subset \R^d$ is a bounded Lipschitz domain and suppose that the family of domain perturbations $\{ \TT\}$ is twice differentiable. 
 
Put $D=\Omega_t$, and let $\Lambda^p(D)$, $1\leq p\leq d$, be the set of $p$-forms on $D$. The outer derivative and wedge product are denoted by $d$ and $\wedge$, respectively. Given 
\[ \alpha=\sum_i \alpha^i dx^i, \ \beta=\sum_i \beta^idx^i\in \Lambda^1(D), \quad x=(x^i), \] 
let  
\[ (\alpha, \beta)=\sum_i \alpha^i\beta^i. \] 
Given  
\[ \lambda=\alpha^1\wedge \cdots \wedge \alpha^p, \ \mu=\beta^1\wedge \cdots \wedge \beta^p \in L^p(D), \] we put 
\[  (\lambda, \mu)=\mbox{det} \ ((\alpha^i, \beta^j))_{ij},  \] 
which is independent of the choice of $\alpha^i, \beta^i\in\Lambda^1(D)$, $1\leq i\leq p$, to represent $\alpha, \beta \in \Lambda^p(D)$. Then the Hodge operator $\ast:\Lambda^p(D)\rightarrow \Lambda^{d-p}(D)$ is defined by 
\[ \omega\wedge \tau=(\ast\omega, \tau)dx^1\wedge \cdots \wedge dx^d, \quad \omega\in \Lambda^p(D), \ \tau\in \Lambda^{d-p}(D). \] 
It holds that 
\[ \ast(dx^{j_1}\wedge \cdots \wedge dx^{j_p})=\mbox{sgn} \ \sigma\cdot dx^{j_{p+1}}\wedge \cdots \wedge dx^{j_d}, \] 
where $\sigma:(1,\cdots, d)\mapsto (j_1, \cdots, j_d)$. By this definition, there arises \cite{flanders} that 
\[ \bnormal ds=(\ast dx_1, \cdots, \ast dx_n), \] 
where $\bnormal=(\nu^i)$ and $ds$ denote the outer unit normal and area element on $\partial D$, respectively. Then we obtain 
\[   
\int_{\partial D}C\cdot \bnormal \ ds=\sum_i\int_{\partial D}C^i\ast dx^i 
\] 
for $C=(C^i)$. 

Let $y=\TT x$ be the transformation of variables, and $\bfa=(a^i) \in C^{1,1}(\Gamma)$ for $\Gamma\subset \R^{d+1}$ defined by (\ref{gamma}). Then we obtain 
\begin{equation} 
\int_{\partial \Omega_t}\bnormal\cdot \bfa \ ds_t=\sum_i\int_{\partial \Omega_t}a^i\ast dy^i. 
 \label{eqna1}
\end{equation} 
In (\ref{T-taylor}) we have 
\[ y^i=x^i+tS^i(x)+\frac{t^2}{2}R^i(x)+o(t^2), \quad  t\rightarrow 0, \ 1\leq i\leq d \] 
in $C^{0,1}(\overline{\Omega})$ and hence 
\[ y^i_j=\delta_{ij}+tS^i_j+\frac{t^2}{2}R^i_j+o(t^2), \quad t\rightarrow 0, \ 1\leq i, j\leq d \] 
in $L^\infty(\Omega)$, where $y^i_j=\frac{\partial y^i}{\partial x^j}$, $S^i_j=\frac{\partial S^i}{\partial x^j}$, $R^i_j=\frac{\partial R^i}{\partial x^j}$, and so forth. Thus it holds that 
\[ dy^i=\sum_jy^i_jdx^j=\sum_j \left(\delta_{ij}+tS^i_j+\frac{t^2}{2}R^i_j\right)dx^j+o(t^2). \]

Pull back the vector field $\bfa(y,t)$ on $\partial\Omega_t$ to that on $\partial\Omega$ by $\TT$: $\bfa(\TT(x),t)$. Using  
\[ 
    a_t^i(y,t) = \frac{\partial a^i}{\partial t}(y,t), \quad
   a_{tt}^i(y,t) = \frac{\partial^2 a^i}{\partial t^2}(y,t), 
\] 
we obtain 
\begin{eqnarray*}
    \frac{\partial \;}{\partial t} a^i(\TT(x),t)
  & = & \nabla_x a^i(\TT(x),t) \frac{\partial\TT(x)}{\partial t}
   + a_t^i(\TT(x),t),  \\ 
   \frac{\partial^2 \;}{\partial t^2} a^i(\TT(x),t)
  & = & \nabla_x^2 a^i(\TT(x),t) \left[
  \frac{\partial\TT(x)}{\partial t},
  \frac{\partial\TT(x)}{\partial t}\right]
  + 2 \nabla_x a_t^i(\TT(x),t) \frac{\partial\TT(x)}{\partial t} \\
  & & + \nabla_x a^i(\TT(x),t)
   \frac{\partial^2\TT(x)}{\partial t^2}
   + a_{tt}^i(\TT(x),t), 
\end{eqnarray*}
and therefore,  
\begin{eqnarray}
a^i(\TT(x),t) & = & a_0^i + t(\Sv\cdot\nabla a_0^i + \dot{a}_0^i)+ \frac{t^2}{2}( \Tv\cdot\nabla a_0^i+ (\nabla^2 a_0^i)[\Sv,\Sv] \nonumber\\ 
& & + 2 \nabla \dot{a}_0^i \Sv + \ddot{a}_0^i)+ o(t^2), 
 \label{72}
\end{eqnarray}
where 
\[  
   a_0^i = a^i(\cdot,0), \quad
   \dot{a}_0^i = a_t^i(\cdot,0), \quad
   \ddot{a}_0^i = a_{tt}^i(\cdot,0). 
\] 

We have, on the other hand,  
\[ 
   *\dd z^i = (-1)^{i-1}\dd z^{1} \wedge \cdots \wedge
      \widehat{\dd z^i} \wedge \cdots \wedge \dd z^{d} 
\] 
for $z = y$ or $z = x$, recalling that $\widehat{\dd z^i}$ indicates the exclusion of $\dd z^i$. It holds also that 
\[ 
   \dd y^i = \sum_{j=1}^d y_j^i \dd x^j, 
\] 
and therefore, 
\begin{eqnarray*} 
*\dd y^i & = & \dd y^{1} \wedge \cdots \wedge
      \widehat{\dd y^i} \wedge \cdots \wedge \dd y^{d} \\ 
& = & (-1)^{i-1}  \sum_{j=1}^d   
     \det (y_p^q)_{\substack{p=1,\cdots,d, p \neq j \\ q = 1,\cdots,d, q \neq i}}
   \dd x^{1} \wedge \cdots \wedge
      \widehat{\dd x^j} \wedge \cdots \wedge \dd x^{d} \\ 
   & = & (-1)^{i-1}  \sum_{j=1}^d (-1)^{j-1}
     \det (y_p^q )_{\substack{p=1,\cdots,d, p \neq j \\ q = 1,\cdots,d, q \neq i}}
     *\!\dd x^j. 
\end{eqnarray*}

First, if $i = j$, we have 
\[ 
   \det (y_p^q)_{\substack{p, q=1,\cdots,d \\
       p, q \neq i}}
  = \det \begin{pmatrix}
       1 + t S_1^1 + \frac{t^2}{2}R_1^1 & 
     t S_2^1 + \frac{t^2}{2}R_2^1 & \cdots &
       t S_d^1 + \frac{t^2}{2}R_d^1  \\
            t S_1^2 + \frac{t^2}{2}R_1^2 & 
      1 + t S_2^2 + \frac{t^2}{2}R_2^2 & \cdots &\
       t S_d^2 + \frac{t^2}{2}R_d^2 \\
    \vdots & \vdots & & \vdots \\
              t S_1^d + \frac{t^2}{2}R_1^d & 
       t S_2^d + \frac{t^2}{2}R_2^d & \cdots &
       1 + t S_d^d + \frac{t^2}{2}R_d^d
    \end{pmatrix}
    + o(t^2),  
\] 
where the entries in the form of $tS_i^q + \frac{t^2}{2}R_i^q$, $tS_p^i + \frac{t^2}{2}R_p^i$, or $1 + tS_i^i + \frac{t^2}{2}R_i^i$ are not included in this matrix. Hence it follows that 
\[ 
\det (y_p^q)_{\substack{p, q=1,\cdots,d \\
       p, q \neq i}}
= 1 + t \sum_{k \neq i} S_k^k + \frac{t^2}{2} \sum_{k \neq i} R_k^k 
    + t^2 \sum_{\substack{j, k \neq i \\ j < k}}
      (S_j^j S_k^k - S_k^j S_j^k) + o(t^2). 
\] 
Second, if  $i \neq j$ we have 
\begin{equation} 
  \det (y_p^q)_{\substack{p=1,\cdots,d, p \neq j
         \\ q=1,\cdots,d, q \neq i}} = (-1)^{j-i+1}
    \left(t S_i^j + \frac{t^2}{2}R_i^j
   + t^2 \sum_{p \neq i,j} (S_i^j S_p^p - S_i^p S_p^j)\right) 
   + o(t^2),    
   \label{to-be-shown}
\end{equation} 
 or 
\[ 
 (-1)^{i+j} \det (y_p^q 
   )_{\substack{p=1,\cdots,d, p \neq j \\ q = 1,\cdots,d, q \neq i}}
  = - t S_{i}^{j} - \frac{t^2}{2} R_{i}^{j}
    + t^2 \sum_{p \neq i,j}
      (S_i^p S_p^j - S_p^p S_i^j) + o(t^2). 
\] 
Equality (\ref{to-be-shown}) is obtained by an expansion of the determinant, and the proof is left for the reader. We thus end up with 
\begin{eqnarray}
  *\dd y^i & = & *\dd x^i + t \sum_{j \neq i} S_j^j *\!\dd x^i
     - t \sum_{j \neq i} S_i^j *\!\dd x^j
     + \frac{t^2}{2} \left(\sum_{j \neq i} R_j^j \right) *\!\dd x^i
     - \frac{t^2}{2} \sum_{j \neq i} R_i^j *\!\dd x^j \nonumber\\
   & & + t^2 \sum_{\substack{j, k \neq i \\ j < k}}
      (S_j^j S_k^k - S_k^j S_j^k) *\!\dd x^i
   + t^2 \sum_{j \neq i}\sum_{p \neq i,j}
      (S_i^p S_p^j - S_p^p S_i^j) *\!\dd x^j + o(t^2). \label{74}
\end{eqnarray}

Writing (\ref{72}) as  
\[ 
    a^i(y,t) = a_0^i + t X^i + \frac{t^2}{2} Y^i + o(t^2), 
\] 
with  
\[ X^i = \Sv\cdot\nabla a_0^i + \dot{a}_0^i, \quad  
Y^i= \Tv\cdot\nabla a_0^i + (\nabla^2 a_0^i)[\Sv,\Sv]
     + 2 \nabla \dot{a}_0^i \cdot \Sv + \ddot{a}_0^i, \]  
we obtain 
\begin{eqnarray*}
   a^i(y,t) *\!\dd y^i & = & a_0^i *\!\dd x^i + 
    t \left(a_0^i \sum_{j \neq i} S_j^j + 
       \Sv\cdot \nabla a_0^i + \dot{a}_0^i\right) *\!\dd x^i
   - t a_0^i \sum_{j \neq i} S_i^j *\!\dd x^j \\
  & & + \frac{t^2}{2} a_0^i \left(\sum_{j \neq i} R_j^j\right) *\!\dd x^i
   - \frac{t^2}{2} a_0^i \sum_{j \neq i} R_i^j *\!\dd x^j \\ 
  & & + t^2 a_0^i \sum_{\substack{j, k \neq i \\ j < k}}
      (S_j^j S_k^k - S_k^j S_j^k) *\!\dd x^i
   + t^2 a_0^i \sum_{j \neq i}\sum_{p \neq i,j}
      (S_i^p S_p^j - S_p^p S_i^j) *\!\dd x^j \\
  & & + \frac{t^2}{2} Y^i *\!\dd x^i
   + t^2 X^i \left(\sum_{j \neq i} S_j^j\right) *\!\dd x^i
  - t^2 X^i \sum_{j \neq i} S_i^j *\!\dd x^j
    + o(t^2)
\end{eqnarray*}
by (\ref{74}). Then (\ref{eqna1}) implies 
\begin{eqnarray}
   \int_{\partial\Omega_t} \bfa(y,t)\cdot\bnormal \dd s_y & = & 
   \sum_{i = 1}^d \int_{\partial\Omega_t} a^i(y,t) *\!\dd y^i \nonumber\\
  & = & \sum_{i = 1}^d \int_{\partial\Omega} a_0^i *\!\dd x^i
    + t \int_{\partial\Omega} I +  
   + \frac{t^2}{2}\int_{\partial\Omega} II + o(t^2) 
 \label{75} 
\end{eqnarray}
with 
\begin{equation} 
  I = \sum_{i=1}^d \sum_{j \neq i} 
   \left[\left(a_0^i S_j^j - a_0^j S_j^i\right)
    + \Sv\cdot \nabla a_0^i + \dot{a}_0^i\right] *\!\dd x^i
 \label{76}
\end{equation} 
and 
\begin{eqnarray}
  II & = & \sum_{i=1}^d a_0^i \left(\sum_{j \neq i} R_j^j\right) *\!\dd x^i
   - \sum_{i=1}^d a_0^i \sum_{j \neq i} R_i^j *\!\dd x^j \nonumber\\
& & + 2 \sum_{i=1}^d a_0^i \sum_{\substack{j, k \neq i \\ j < k}}
      \left(S_j^j S_k^k - S_k^j S_j^k\right) *\!\dd x^i
   + 2 \sum_{i=1}^d a_0^i \sum_{j \neq i}\sum_{p \neq i,j}
      \left(S_i^p S_p^j - S_p^p S_i^j\right) *\!\dd x^j \nonumber\\ 
& & + \sum_{i=1}^d
   (\Tv\cdot\nabla a_0^i + (\nabla^2 a_0^i)[\Sv,\Sv]
     + 2 \nabla \dot{a}_0^i \cdot \Sv + \ddot{a}_0^i) *\!\dd x^i \nonumber\\ 
&  & + 2 \sum_{i=1}^d (\Sv\cdot\nabla a_0^i + \dot{a}_0^i)
   \left(\sum_{j \neq i} S_j^j\right) *\!\dd x^i
  - 2 \sum_{i=1}^d (\Sv\cdot\nabla a_0^i + \dot{a}_0^i)
   \sum_{j \neq i} S_i^j *\!\dd x^j. 
   \label{77x} 
\end{eqnarray}

Here, we obtain     
\begin{eqnarray*}
  I & = & \sum_{i=1}^d \left[ \sum_{j = 1}^d 
   (a_0^i S_j^j - a_0^j S_j^i)
    + \sum_{j = 1}^d S^j a_{0,j}^i
    - \sum_{j=1}^d a_{0,j}^j S^i + \sum_{j=1}^d a_{0,j}^j S^i
    \right] *\!\dd x^i
   + \sum_{i = 1}^d \dot{a}_0^i *\!\dd x^i \\ 
  & = & \sum_{i=1}^d \left[ \sum_{j = 1}^d 
   (a_0^i S^j)_j - 
    \sum_{j=1}^d (a_0^j S^i)_j
   + \sum_{j=1}^d a_{0,j}^j S^i \right] *\!\dd x^i
   + \sum_{i = 1}^d \dot{a}_0^i *\!\dd x^i \\
   & = &   \sum_{i=1}^d \left[(\nabla\cdot\bfa)S^i +
     \dot{a}_0^i \right] *\!\dd x^i
      + \theta_1
\end{eqnarray*}
for 
\[ 
  \theta_1= \sum_{i = 1}^d \left[ \sum_{j = 1}^d 
   (a_0^i S^j)_j - 
   \sum_{j=1}^d (a_0^j S^i)_j\right] *\!\dd x^i. 
\] 
Since $dx^i\wedge \ast dx^i=dx^1\wedge \cdots \wedge dx^d$, it holds that  
\[ 
  \dd \theta_1 = \sum_{i,j = 1}^d \left[ 
   (a_0^i S^j)_{ij} - 
    (a_0^j S^i)_{ij}\right] 
   \dd x^1 \wedge \cdots \wedge \dd x^d  = 0, 
\] 
and hence 
\[ 
   \int_{\partial\Omega} \theta_1 = \int_\Omega \dd \theta_1 = 0
\] 
by the Stokes theorem. We thus obtain 
\begin{equation}
  \int_{\partial\Omega} I = 
 \sum_{i=1}^d \int_{\partial\Omega} \left[(\nabla\cdot\bfa)S^i +
     \dot{a}_0^i \right] *\!\dd x^i 
  =  \int_{\partial\Omega} \left[(\nabla\cdot\bfa)\Sv \cdot\bnormal
   +   \dot{\bfa}_0\cdot\bnormal \right] \dd s.  
 \label{78}
\end{equation}

We divide $II$ into four terms, involving $\ddot{\bfa}_0$, $\{\Tv, \bfa_0\}$, $\{\Sv, \dot{\bfa}_0\}$, and $\{ \Sv, \bfa_0\}$, denoted by $II_1$, $II_2$, $II_3$, and $II_4$, respectively. First, we have  
\[ 
  II_1 = \sum_{i=1}^d \ddot{a}_0^i *\!\dd x^i, 
\] 
and hence 
\[ 
\int_{\partial\Omega} II_1 
   = \int_{\partial\Omega} \ddot{\bfa}_0 \cdot \bnormal \dd s. 
\] 
Second, there arises that 
\begin{eqnarray*} 
  II_2 & = & \sum_{i=1}^d a_0^i 
    \left(\sum_{j \neq i} R_j^j\right) *\!\dd x^i
   - \sum_{i=1}^d a_0^i \sum_{j \neq i} R_i^j *\!\dd x^j
     + \sum_{i=1}^d \Tv\cdot\nabla a_0^i *\!\dd x^j \\ 
 & = & \sum_{i=1}^d \sum_{j \neq i} 
   \left[(a_0^i R_j^j - a_0^j R_j^i)
    + \Tv\cdot \nabla a_0^i\right] *\!\dd x^i \\ 
   & = & \sum_{i=1}^d \left[ \sum_{j = 1}^d 
   (a_0^i R_j^j - a_0^j R_j^i)
    + \sum_{j = 1}^d R^j a_{0,j}^i
    - \sum_{j=1}^d a_{0,j}^j R^i + \sum_{j=1}^d a_{0,j}^j R^i
    \right] *\!\dd x^i \\ 
   & = & \sum_{i=1}^d \left[ \sum_{j = 1}^d 
   (a_0^i R^j)_j - 
    \sum_{j=1}^d(a_0^j R^i)_j
   + \sum_{j=1}^d a_{0,j}^j R^i \right] *\!\dd x^i \\
   & = & \sum_{i=1}^d (\nabla\cdot\bfa)R^i *\!\dd x^i + \theta_2 
\end{eqnarray*}
with 
\[ 
  \theta_2= \sum_{i = 1}^d \left[ \sum_{j = 1}^d 
   (a_0^i R^j)_j - 
   \sum_{j=1}^d (a_0^j R^i)_j\right] *\!\dd x^i, 
\] 
Then we obtain $\dd \theta_2 = 0$, and hence 
\[ 
   \int_{\partial\Omega} \theta_2 = 0,
\] 
similarly. It thus follows that 
\[ 
  \int_{\partial\Omega} II_2 =
   \int_{\partial\Omega} 
  (\nabla\cdot\bfa)\left(\Tv\cdot \bnormal\right) \dd s. 
\] 
Third, we have 
\begin{eqnarray*}
  II_3 & = & 2 \sum_{i=1}^d \dot{a}_0^i 
    \left(\sum_{j \neq i} S_j^j \right) *\!\dd x^i
   - 2 \sum_{i=1}^d \dot{a}_0^i \sum_{j \neq i} S_i^j *\!\dd x^j
     + 2 \sum_{i=1}^d \Sv\cdot\nabla \dot{a}_0^i *\!\dd x^j \\
 & = & 2 \sum_{i=1}^d \sum_{j \neq i} 
   \left[(\dot{a}_0^i S_j^j - \dot{a}_0^j S_j^i)
    + \Sv\cdot \nabla \dot{a}_0^i\right] *\!\dd x^i \\ 
   & = & 2 \sum_{i=1}^d \left[ \sum_{j = 1}^d 
   (\dot{a}_0^i S_j^j - \dot{a}_0^j S_j^i)
    + \sum_{j = 1}^d S^j \dot{a}_{0,j}^i
    - \sum_{j=1}^d \dot{a}_{0,j}^j S^i
    + \sum_{j=1}^d \dot{a}_{0,j}^j S^i
    \right] *\!\dd x^i \\ 
   & = & 2 \sum_{i=1}^d \left[ \sum_{j = 1}^d 
   (\dot{a}_0^i S^j)_j - 
    \sum_{j=1}^d (\dot{a}_0^j S^i)_j
   + \sum_{j=1}^d \dot{a}_{0,j}^j S^i \right] *\!\dd x^i \\ 
   & = & 2 \sum_{i=1}^d (\nabla\cdot\dot{\bfa})S^i *\!\dd x^i + 2 \theta_3  
\end{eqnarray*}
for 
\[ 
  \theta_3= \sum_{i = 1}^d \left[ \sum_{j = 1}^d 
   (\dot{a}_0^i S^j)_j - 
   \sum_{j=1}^d(\dot{a}_0^j S^i)_j\right] *\!\dd x^i. 
\] 
Then we obtain $\dd \theta_3 = 0$, and hence 
\[ 
\int_{\partial\Omega} \theta_3 = 0, 
\] 
similarly. It thus follows that 
\[ 
\int_{\partial\Omega} II_3 =
   2 \int_{\partial\Omega} 
  (\nabla\cdot\dot{\bfa})\left(\Sv\cdot \bnormal\right) \dd s. 
\] 

Finally, we treat 
\begin{eqnarray*}
& & II_4 =  2 \sum_{i=1}^d a_0^i \sum_{\substack{j, k \neq i \\ j < k}}
      (S_j^j S_k^k - S_k^j S_j^k) *\!\dd x^i
   + 2 \sum_{i=1}^d a_0^i \sum_{j \neq i}\sum_{p \neq i,j}
      (S_i^p S_p^j - S_p^p S_i^j) *\!\dd x^j \\ 
   & & + \sum_{i=1}^d
   (\nabla^2 a_0^i)[\Sv,\Sv] *\!\dd x^i 
   + 2 \sum_{i=1}^d (\Sv\cdot\nabla a_0^i)
   \left(\sum_{j \neq i} S_j^j\right) *\!\dd x^i
  - 2 \sum_{i=1}^d (\Sv\cdot\nabla a_0^i)
   \sum_{j \neq i} S_i^j *\!\dd x^j 
\end{eqnarray*}
Our goal is to show 
\begin{eqnarray}
  \int_{\partial\Omega} II_{4} & = & \sum_{i=1}^{d}
  \int_{\partial\Omega} \nabla\cdot[(\nabla\cdot\bfa_{0})\Sv]
  S^{i} *\!\dd x^{i} - \sum_{i=1}^{d} \int_{\partial\Omega}
  (\nabla\cdot\bfa_{0}) (S\cdot\nabla)S^{i} *\!\dd x^{i} \nonumber\\
  & = & 
  \int_{\partial\Omega} \nabla\cdot[(\nabla\cdot\bfa)\Sv]
  (\Sv\cdot\bnormal)  - 
  (\nabla\cdot\bfa) [(S\cdot\nabla) \Sv]\cdot \bnormal \ \dd s. 
   \label{79}
\end{eqnarray}
For this purpose we note the equalities 
\begin{eqnarray*}
& & \nabla\cdot [(\nabla\cdot\bfa_{0})\Sv]
  S^{i} *\!\dd x^{i} - 
    (\nabla\cdot\bfa_{0}) (S\cdot\nabla)S^{i} *\!\dd x^{i} \\ 
& & \quad = \sum_{j,k} a_{0,jk}^{k} S^{j}S^{i} *\!\dd x^{i}
     + \sum_{j,k} a_{0,k}^{k} S_{j}^{j}S^{i} *\!\dd x^{i}
     - \sum_{j,k} a_{0,k}^{k} S_{j}^{i}S^{j} *\!\dd x^{i} 
\end{eqnarray*}
and 
\[ 
  \nabla^{2} a_{0}^{i} [\Sv,\Sv]
      = \sum_{j,k} a_{0,jk}^{i} S^{j}S^{k}, \quad
        \Sv\cdot\nabla a_{0}^{i} = \sum_{k} a_{0,k}^{i}S^{k}. 
\] 
We write also  
\begin{eqnarray*}
& &  2 \sum_{i} a_0^i \sum_{\substack{j, k \neq i \\ j < k}}
      (S_j^j S_k^k - S_k^j S_j^k) *\!\dd x^i
      = 2 \sum_{i} a_{0}^{i} A^{i} *\!\dd x^{i} , \\      
& & 2 \sum_{i} a_0^i \sum_{j \neq i}\sum_{p \neq i,j}
      (S_i^p S_p^j - S_p^p S_i^j) *\!\dd x^j
      = 2 \sum_{i} a_{0}^{i} \sum_{j \neq i} B^{ij} *\!\dd x^{j}, 
\end{eqnarray*} 
using 
\begin{eqnarray*}
A^i & = & \sum_{\substack{j, k \neq i \\ j < k}}
( S_j^j S_k^k - S_k^j S_j^k)
= \frac{1}{2}\sum_{j \neq i} \sum_{k \neq i}
(S_j^j S_k^k - S_k^j S_j^k), \\
B^{ij} & = & \sum_{k \neq i,j}
(S_i^k S_k^j - S_k^k S_i^j)
= \sum_{k \neq i}(S_i^k S_k^j - S_k^k S_i^j). 
\end{eqnarray*}
We thus obtain 
\[ 
  II_{4} = \sum_{i} \nabla\cdot[(\nabla\cdot\bfa_{0})\Sv]
  S^{i} *\!\dd x^{i} - \sum_{i}
    (\nabla\cdot\bfa_{0}) (S\cdot\nabla)S^{i} *\!\dd x^{i}
   + \theta_{4} + \theta_{5} + \theta_{6},
\] 
for 
\[ 
 \theta_{4} = \sum_{i} \sum_{j,k}
     a_{0,jk}^{i}S^{j}S^{k} *\!\dd x^{i}
     - \sum_{i} \sum_{j,k}
       a_{0,jk}^{k}S^{i}S^{j} *\!\dd x^{i}
     - \sum_{i} \sum_{j,k} a_{0,k}^{k} S^{i}S_{j}^{j}
        *\!\dd x^{i}, 
\] 
\begin{eqnarray*}
  \theta_{5} & = & \sum_{i} \sum_{j,k} a_{0,k}^{k} S^{j} S_{j}^{i}
             *\!\dd x^{i} 
            + 2 \sum_{i} \sum_{j,k, j \neq i} 
             a_{0,k}^{i}S^{k}S_{j}^{j} *\!\dd x^{i}
             - 2 \sum_{i} \sum_{j,k,j\neq i}
              a_{0,k}^{i}S^{k} S_{i}^{j} *\!\dd x^{j} \\
      & = & \sum_{i} \sum_{j,k} a_{0,k}^{k} S^{j} S_{j}^{i}
             *\!\dd x^{i} 
           + 2 \sum_{i} \sum_{j,k} 
             a_{0,k}^{i}S^{k}S_{j}^{j} *\!\dd x^{i}
             - 2 \sum_{i} \sum_{j,k}
              a_{0,k}^{j}S^{k} S_{j}^{i} *\!\dd x^{i}, 
\end{eqnarray*} 
and 
\[ 
   \theta_{6} = 2 \sum_{i} a_{0}^{i} A^{i} *\!\dd x^{i} +
         2 \sum_{i} a_{0}^{i} \sum_{j \neq i} B^{ij} *\!\dd x^{j}. 
\] 
Thus, equality (\ref{79}) is reduced to 
\[ 
  \int_{\partial\Omega} (\theta_{4} + \theta_{5} + \theta_{6}) 
   = \int_{\Omega} \dd(\theta_{4} + \theta_{5} + \theta_{6}) = 0, 
\] 
which follows from 
\begin{equation} 
\dd(\theta_{4}+\theta_{5}+\theta_{6})=0. 
 \label{80}
\end{equation} 

In fact, we have 
\[ \dd (\theta_4+\theta_5+\theta_6)=X \ \dd x^1 \wedge \cdots \wedge \dd x^d \] 
with 
\begin{eqnarray*}
X & = & \sum_{i,j,k} a_{0,ijk}^{i} S^{j}S^{k}
  + \sum_{i,j,k} a_{0,jk}^{i} (S^{j}S^{k})_{i}
  - \sum_{i,j,k} a_{0,ijk}^{k} S^{i}S^{j}
  - \sum_{i,j,k} a_{0,jk}^{k}(S^{i}S^{j})_{i} \\
 & &  - \sum_{i,j,k} a_{0,ik}^{k} S^{i}S_{j}^{j} 
    - \sum_{i,j,k} a_{0,k}^{k} (S^{i}S_{j}^{j})_{i} \\
& & + \sum_{i,j,k} a_{0,ik}^{k} S^{j}S_{j}^{i}
  + \sum_{i,j,k} a_{0,k}^{k} (S^{j}S_{j}^{i})_{i}
  + 2 \sum_{i,j,k} a_{0,ik}^{i} S^{k}S_{j}^{j}
  + 2 \sum_{i,j,k} a_{0,k}^{i} (S^{k}S_{j}^{j})_{i} \\
& &   - 2 \sum_{i,j,k} a_{0,ik}^{j} S^{k}S_{j}^{i}
  - 2 \sum_{i,j,k} a_{0,k}^{j} (S^{k}S_{j}^{i})_{i} \\
& & + 2\sum_{i} a_{0,i}^{i}A_{i} 
  + 2 \sum_{i} a_{0}^{i} A_{i}^{i}
  + 2 \sum_{i} \sum_{j \neq i} a_{0,j}^{i}B^{ij}
  + 2 \sum_{i} a_{0}^{i} \sum_{j \neq i} B_{j}^{ij} \\ 
  & \equiv & X_1+X_2-X_3-X_4-X_5-X_6 \\ 
  & & +X_7+X_8+2X_{9}+2X_{10}-2X_{11}-2X_{12} \\ 
  & & +2X_{13}+2X_{14}+2X_{15}+2X_{16}. 
\end{eqnarray*}
First, noticing the terms involving the third order derivatives of $a_{0}^{i}$, we realize 
\[ X_1-X_3=0. \] 
Second, the terms involving the zero-th order derivatives of $a_0^i$ also cancel as 
\[ X_{14}+X_{16}=\sum_ia_0^iF_i \] 
with 
\begin{eqnarray*}
F_i & = & A_{i}^{i} + \sum_{j \neq i} B_{j}^{ij} \\ 
& = & \frac{1}{2} \sum_{j, k \neq i}
   (S_j^j S_k^k - S_k^j S_j^k)_i
    + \sum_{j,k \neq i} (S_i^k S_k^j - S_k^k S_i^j)_j \\
  & = & \frac{1}{2}\sum_{j,k \neq i}
   ( S_{ij}^j S_k^k + S_j^j S_{ki}^k 
    - S_{ik}^j S_j^k - S_k^j S_{ij}^k)
    + \sum_{j, k \neq i}
      (S_{ij}^k S_k^j
    + S_i^k S_{jk}^j
    - S_{kj}^k S_i^j 
    - S_k^k S_{ij}^j) \\ 
   & = & \frac{1}{2}(X_{17}+X_{18}-X_{19}-X_{20})+X_{21}+X_{22}-X_{23}-X_{24}.  
\end{eqnarray*} 
Since 
\[ X_{17}=X_{24}, \quad X_{20}=X_{21}, \quad X_{22}=X_{23}, \] 
we obtain 
\[  
 F_i = \frac{1}{2}\sum_{j,k \neq i} 
   (- S_{ij}^j S_k^k + S_j^j S_{ki}^k 
    - S_{ik}^j S_j^k + S_k^j S_{ij}^k) = 0. \] 
Third, the terms involving the second order derivatives of $a_{0}^{i}$ cancel as
\begin{eqnarray*}
& & X_2-X_4-X_5+X_7+2X_9-2X_{11}  \\ 
& & \quad = \sum_{i,j,k}a_{0,jk}^{i} S_{i}^{j}S^{k}
    + \sum_{i,j,k} a_{0,jk}^{i} S^{j}S_{i}^{k}
  + \sum_{i,j,k} a_{0,ik}^{k} S^{j}S_{j}^{i}
  + 2 \sum_{i,j,k} a_{0,ik}^{i} S^{k}S_{j}^{j} \\ 
& & \qquad  - \sum_{i,j,k} a_{0,jk}^{k} S^{j}S_{i}^{i}
   - \sum_{i,j,k} a_{0,jk}^{k} S^{i}S_{i}^{j}
   - \sum_{i,j,k} a_{0,ik}^{k} S^{i}S_{j}^{j}
   - 2 \sum_{i,j,k} a_{0,ik}^{j} S^{k}S_{j}^{i} \\ 
   & & \quad = X_{25}+X_{26}+X_{27}+2X_{28}-X_{29}-X_{30}-X_{31}-2X_{32} =0 
\end{eqnarray*}
by 
\[ X_{25}=X_{26}=X_{32}, \quad X_{27}=X_{30}, \quad X_{28}=X_{29}=X_{31}. \] 
Finally, for the terms involving the first order derivatives of $a_0^i$ we obtain  
\begin{eqnarray*}
& &Y= -X_6+X_8+2X_{10}-2X_{12}+2X_{13} \\ 
& & \quad =  - \sum_{i,j,k} a_{0,k}^{k} (S^{i}S_{j}^{j})_{i} 
  + \sum_{i,j,k} a_{0,k}^{k} (S^{j}S_{j}^{i})_i
  + 2 \sum_{i,j,k} a_{0,k}^{i} (S^{k}S_{j}^{j})_{i} \\
 & & \qquad - 2 \sum_{i,j,k} a_{0,k}^{j} (S^{k}S_{j}^{i})_{i}
  + 2 \sum_{i} a_{0,i}^{i}A_{i} 
  + 2 \sum_{i} a_{0,j}^{i} \sum_{j \neq i}B^{ij} \\
 & & \quad = - \sum_{i,j,k} a_{0,k}^{k} S_{i}^{i}S_{j}^{j}
     - \sum_{i,j,k} a_{0,k}^{k} S^{i}S_{ij}^{j}
    + \sum_{i,j,k} a_{0,k}^{k} S_{i}^{j}S_{j}^{i}
    + \sum_{i,j,k} a_{0,k}^{k} S^{j}S_{ij}^{i} \\
&  & \qquad + 2
        \sum_{i,j,k} a_{0,k}^{i} S_{i}^{k}S_{j}^{j}
    + 2 \sum_{i,j,k} a_{0,k}^{i} S^{k}S_{ij}^{j}
    - 2 \sum_{i,j,k} a_{0,k}^{j} S_{i}^{k}S_{j}^{i}
    - 2 \sum_{i,j,k} a_{0,k}^{j} S^{k}S_{ij}^{i} \\
&  & \qquad + 2 \sum_{i} a_{0,i}^{i}A_{i} 
  + 2 \sum_{i} \sum_{j \neq i} a_{0,j}^{i} B^{ij} \\ 
& & =-X_{33}-X_{34}+X_{35}+X_{36} \\ 
& & \quad +2X_{37}+2X_{38}-2X_{39}-2X_{40}+2X_{41}+2X_{42}. 
\end{eqnarray*} 
Then it holds that 
\begin{eqnarray*} 
& & -X_{33}-X_{34}+X_{35}+X_{36}  = \sum_{k} a_{0,k}^{k} \left(\sum_{i,j \neq k}
        (S_{i}^{i}S_{j}^{j} - S_{i}^{j}S_{j}^{i}) 
    + \sum_{i,j} S_{i}^{j}S_{j}^{i}
   - \sum_{i,j} S_{i}^{i}S_{j}^{j}\right) \\ 
   & & \quad = 2 \sum_{j,k} a_{0,k}^{k}
         (S_{k}^{j}S_{j}^{k} - S_{k}^{k}S_{j}^{j}), \\ 
& & 2X_{37}-2X_{39}=2 \sum_{i,j,k} a_{0,j}^{i} 
      (S_{i}^{j}S_{k}^{k} - S_{k}^{j}S_{i}^{k})+ 2 \sum_{i} \sum_{j \neq i} a_{0,j}^{i}\sum_{k} (S_{i}^{k}S_{k}^{j} - S_{k}^{k}S_{i}^{j}) \\ 
& & 2X_{38}-2X_{40}= 2 \sum_{i} \sum_{j \neq i} a_{0,j}^{i}\sum_{k}(S_{i}^{k}S_{k}^{j} - S_{k}^{k}S_{i}^{j}) \\ 
& & 2X_{41}+2X_{42}=0, 
\end{eqnarray*} 
which implies 
\[ Y= 2 \sum_{j,k} a_{0,k}^{k}
         (S_{k}^{j}S_{j}^{k} - S_{k}^{k}S_{j}^{j}) 
    + 2 \sum_{i,k} a_{0,i}^{i} 
      (S_{i}^{i}S_{k}^{k} - S_{k}^{i}S_{i}^{k}) = 0. \] 
Hence we obtain (\ref{80}). 

We thus end up with 
\begin{eqnarray*}
\int_{\partial\Omega_t} \bfa(y,t)\cdot\bnormal \dd s_y  & =  & 
\sum_{i = 1}^d \int_{\partial\Omega_t} a^i(y,t) *\!\dd y^i \\ 
& = &  \int_{\partial\Omega} \bfa_0 \cdot \bnormal \dd s
   + t \int_{\partial\Omega} [(\nabla\cdot\bfa)\Sv \cdot\bnormal
   +   \dot{\bfa}_0\cdot\bnormal] \dd s \\ 
& & + 
  \frac{t^2}{2}\int_{\partial\Omega} 
 [(\nabla\cdot
    [(\nabla\cdot\bfa_0)\Sv] + 2 \nabla \cdot \dot{\bfa}_0)\Sv\cdot\bnormal + \ddot{\bfa}_0 \cdot \bnormal] \\ 
& & \quad + 
(\Tv - (\Sv\cdot\nabla)\Sv)\cdot\bnormal \ \dd s + o(t^2),  
\end{eqnarray*}
and hence reach the formulae in Lemmas \ref{vector-first} and \ref{vector-second}.

\section{Second fundamental form on $\partial\Omega$}\label{seca3}

In (\ref{59}) if $\{ \TT\}$ is the normal perturbation (\ref{nperturbation}) we have 
\[ (\nabla \nu)[S,S]=\delta\rho \frac{\partial \nu}{\partial \nu}=0 \] 
by (\ref{consistent}) and hence 
\[ \chi =-(\delta\rho)^2\nabla\cdot \nu, \quad \sigma=0. \] 
If $\{ \TT\}$ is the dynamical perturbation using (\ref{7}), it follows that 
\[ \chi=-(\delta\rho)^2\nabla\cdot\nu + \left(v\cdot \nu\frac{\partial}{\partial\nu}-v\cdot\nabla_\tau\right)\delta\rho, \quad \sigma=\left(v\cdot\nu\frac{\partial}{\partial\nu}-v\cdot\nabla_\tau\right)\delta\rho \] 
from (\ref{xx3}), (\ref{xx4}), and (\ref{19}).  
In \cite{gs52} the second Hadamard variation for $d=3$ is described in accordance with the second fundamental form of $\partial\Omega$. It is concerned on the normal perturbation (\ref{nperturbation}), where $\rho=S\cdot \bnormal$ and $R=0$ 

These values $\chi, \sigma$ used in Theorem \ref{thm21} are actually associated with the second fundamental form  on $\partial\Omega$, defined by 
\[ {\cal B}[\xi,\eta]=-\frac{\partial \nu}{\partial \xi}\cdot \eta=-(\nabla \bnormal)[\xi,\eta], \quad \xi, \eta \in \R^d. \] 
The argument relies on the following formula. Recall 
\[ 
\delta \rho=S\cdot\bnormal, \quad  
\mbox{tr} \ {\cal B}=-\nabla\cdot \bnormal, \quad S_\tau=(S\cdot\tau)\tau. 
 \]  

\begin{lemma}\cite[(3.1.1.8), p.137]{grisvard}\label{lem22} 
It holds that 
\begin{eqnarray} 
& & (\nabla\cdot S)\delta\rho-[(S\cdot\nabla)S]\cdot \bnormal \nonumber\\ 
& & \quad =\nabla_\tau\cdot (\delta \rho S_\tau) -(\mbox{tr} \ {\cal B})(\delta \rho)^2-{\cal B}(S_\tau, S_\tau)-2(S_\tau\cdot \nabla_\tau)\delta\rho.  
 \label{83}
\end{eqnarray}  
\end{lemma} 

\begin{lemma}\label{lem29}
It holds that 
\begin{eqnarray*} 
& & [(S\cdot \nabla S)]\cdot \nu + \left(-\delta \rho
				    \frac{\partial}{\partial
				    \nu}+S_\tau\cdot \nabla_\tau \right)\delta\rho \\ 
& & \quad = -(\nabla\cdot\nu)(\delta\rho)^2+{\mathcal B}(S_\tau, S_\tau)+2(S_\tau\cdot \nabla_\tau)\delta\rho
\end{eqnarray*} 
\end{lemma}

\begin{proof} 
The result follows from a direct computation. First, Lemma \ref{lem22} implies 
\begin{eqnarray*} 
& & (S_\tau\cdot \nabla_\tau)\delta\rho+{\cal B}(S_\tau, S_\tau)+(\nabla\cdot S)\delta \rho \nonumber\\ 
& & \quad = [(S\cdot\nabla)S]\cdot \bnormal +\nabla_\tau\cdot (\delta\rho S_\tau)+(\nabla\cdot \bnormal)(\delta\rho)^2-(S_\tau\cdot \nabla_\tau)\delta\rho  \\ 
& & \quad =[(S\cdot \nabla) S]\cdot \nu +(\nabla\cdot \nu)(\delta \rho)^2+\delta \rho\nabla_\tau\cdot S_\tau 
 \label{81pre}
\end{eqnarray*}  
and hence 
\begin{eqnarray} 
& &[(S\cdot\nabla) S]\cdot \bnormal +(\nabla\cdot\bnormal)(\delta\rho)^2+(S\cdot\nabla)\delta\rho \nonumber\\ 
& & \quad = -\delta \rho \nabla_\tau\cdot S_\tau+(S_\tau\cdot \nabla_\tau)\delta \rho+ {\mathcal B}(S_\tau, S_\tau)+(\nabla\cdot S)\delta \rho +(S\cdot \nabla)\delta \rho \nonumber\\ 
& &\quad = (S_\tau\cdot\nabla_\tau)\delta\rho+{\mathcal B}(S_\tau, S_\tau)-\delta\rho\nabla_\tau\cdot S_\tau+\nabla\cdot (S\delta\rho). 
 \label{81}
\end{eqnarray}
Then we use 
\[ 
(S\cdot \nabla)\delta \rho = \left(\delta\rho \frac{\partial}{\partial \nu}
  + S_\tau\cdot \nabla_\tau \right)\delta\rho 
\] 
and 
\begin{eqnarray*}  
\nabla \cdot (S\delta \rho) & = & \nabla_\tau(S_\tau\delta\rho)+\frac{\partial}{\partial \nu}(\delta \rho)^2 \\ 
& = & \delta \rho \nabla_\tau\cdot S_\tau+(S_\tau\cdot \nabla_\tau)\delta\rho+2\delta\rho \frac{\partial \delta\rho}{\partial \nu},  
\end{eqnarray*}  
which implies 
\[
  [(S\cdot \nabla)S]\cdot \nu+(\nabla\cdot\nu)(\delta\rho)^2
  = {\mathcal  B}(S_\tau,  S_\tau) + 
   \left(S_\tau\cdot\nabla_\tau+\delta\rho
  \frac{\partial}{\partial \nu} \right)\delta\rho 
\] 
by (\ref{81}). Then the result follows. 
\end{proof}  

\begin{theorem} 
It holds that 
\begin{eqnarray*} 
& &\chi=\delta^2\rho-2(S_\tau\cdot\nabla_\tau)\delta\rho-{\cal B}(S_\tau, S_\tau) -(\delta\rho)^2(\nabla\cdot \bnormal) \\ 
& & \sigma=\delta^2\rho-2(S_\tau\cdot\nabla_\tau)\delta\rho-{\cal B}(S_\tau, S_\tau)
\end{eqnarray*} 
in Theorem \ref{thm21}. 
\end{theorem} 

\begin{proof} 
Equalities (\ref{xx3}) and  (\ref{xx4}) imply   
\[ 
\chi = [R-(S\cdot \nabla)S]\cdot \nu - 
   (\delta\rho)^2\nabla\cdot \nu + 
  \left(\delta\rho \frac{\partial}{\partial \nu} -(S_\tau\cdot\nabla_\tau)\right)\delta \rho 
\] 
and 
\[  
\sigma=[R-(S\cdot \nabla)S]\cdot \nu +
  \left(\delta\rho\frac{\partial}{\partial \bnormal}-(S_\tau\cdot\nabla_\tau)\right)\delta\rho,  
\] 
respectively. Then we obtain the result by (\ref{83}). 
 
\end{proof}


\vspace{2cm} 

{\bf Acknowledgements.} This work was promoted in Research
Institute for Mathematical Sciences (RIMS) Joint Research Program
during 2019--2021.  (RIMS is an International Joint Usage/Research
Center located in Kyoto University.)
The  authors thank Professors Hideyuki Azegami and Erika Ushikoshi for many
detailed discussions at these occasions. The first author was supported
by JSPS Grant-in-Aid for Scientific Research 19H01799.  The second
author was supported by JSPS Grant-in-Aid for Scientific Research 21K03372.

\end{document}